\documentclass[]{interact}
\usepackage{lmodern}
\usepackage{epstopdf}
\usepackage[caption=false]{subfig}
\usepackage{hyperref}
\usepackage[numbers,sort&compress]{natbib}
\bibpunct[, ]{[}{]}{,}{n}{,}{,}
\makeatletter
\def\NAT@def@citea{\def\@citea{\NAT@separator}}
\makeatother

\theoremstyle{plain}
\newtheorem{theorem}{Theorem}[section]
\newtheorem{lemma}[theorem]{Lemma}
\newtheorem{corollary}[theorem]{Corollary}
\newtheorem{proposition}[theorem]{Proposition}

\theoremstyle{definition}
\newtheorem{definition}[theorem]{Definition}

\numberwithin{equation}{section}
\theoremstyle{remark}
\newtheorem{remark}{Remark}

\newcommand{\be}{\begin{equation}}
\newcommand{\ee}{\end{equation}}

\newcommand{\bea}{\begin{eqnarray}}
\newcommand{\eea}{\end{eqnarray}}

\begin{document}

\title{A determinant approach for generalized $q$-Bernoulli polynomials and asymptotic results}

\author{
\name{S.Z.H. Eweis\textsuperscript{a}\thanks{Contact  S.Z.H.  Eweis. Email: Sahar.Zareef@bsu.edu.eg } and  Z.S.I. Mansour\textsuperscript{b}\thanks{Contact Z.S.I.  Mansour. Email:zsmansour@cu.edu.eg}}
\affil{\textsuperscript{a} Mathematics and Computer Science Department,  Faculty of Science, Beni-Suef University, Beni-Suef, Egypt.\\
 \textsuperscript{b}Mathematics Department, Faculty of Science, Cairo University, Giza,  Egypt.
}}

\maketitle
\begin{abstract}
 In earlier work, we introduced three families of polynomials where the generating function of each set includes one of the three Jackson $q$-analogs of the Bessel function. This paper gives determinant representation for  each family, their large $n $ asymptotics, and two expansion theorems for specific classes of entire functions. We include two examples.
\end{abstract}

\begin{amscode}
 11B68, 15A15, 32A15, 30E15
\end{amscode}

\section{Priliminaries}\label{Sec.1}

In \cite{sahar}, Ewies and Mansour introduced three sets of polynomials
\begin{equation*}
  (B^{(k)}_{n,\alpha}(\cdot;q))_{n}\quad(0<q<1,\, \alpha >-1,\, k\in\{1,2,3\}).
\end{equation*}
The generating functions of $(B^{(k)}_{n,\alpha}(\cdot;q))_{n}$  include the Jackson  $q$-Bessel functions $J_{\alpha}^{(k)}(\cdot;q)\,\,(k=1,2,3)$ defined in  (\ref{s765}) below. These polynomials are $q$-analogs of Frappier$^{,}$s generalized polynomials $(B_{n,\alpha}(\cdot))_{ n}$ defined by, see \cite{Frappier1, Frappier2, F3},
\begin{equation}\label{q9321}
   \frac{e^{(z-\frac{1}{2})t}}{g_{\alpha}(\frac{it}{2})}=\sum_{n=0}^{\infty}B_{n,\alpha}(z)\frac{t^{n}}{n!},\quad  |t|<2 j_{1,\alpha},
 \end{equation}
where
\begin{equation}\label{q421}
 g_{\alpha}(t)=2^{\alpha}\Gamma(\alpha+1)\frac{ J_{\alpha}(t)}{t^{\alpha}}.
\end{equation}
Here $ J_{\alpha}(t)$ is the Bessel function of the first kind of order $\alpha,$ and $j_{1,\alpha} $ is the smallest positive zero of $J_{\alpha}(t)$. In \cite{sahar}, we studied the main properties of the polynomials $(B^{(k)}_{n,\alpha}(\cdot;q))_{n}$ and their large $n$-degree asymptotics.  This paper is organized as follows. In the remainder of this section, we introduce the $q$-notations and notions needed in the sequel. In Section 2, we investigate a determinant representation of the generalized $q$-Bernoulli polynomials. Section 3 includes their large $n$-degree asymptotics. Finally in  Section 4, we  expand a special class of entire functions in terms of the generalized $q$-Bernoulli polynomials associated with the second Jackson $q$-Bessel function. Throughout this paper, we use $\mathbb{N}$  to indicate the set of natural numbers, $\mathbb{N}_{0}:=\mathbb{N}\cup\{0\}$ indicates the set of non negative integers and $\mathbb{C}$ denotes the set of complex numbers. We refer the readers to  \cite{Gasper} for definitions of the  $q$-shifted factorial, the $q$-natural number, the $q$-factorial,  the $q$-binomial coefficients, the basic hypergeometric functions and the $q$-gamma function. Jackson in \cite{Jackson1} defined the $q$-difference operator by
       \begin{equation*}
           D_q f(z)= \dfrac{f(qz)-f(z)}{z(q-1)}, \quad z\neq 0.
         \end{equation*}

The symmetric $q$-difference operator is defined by, see \cite{Cardoss, Gasper},
\begin{equation}\label{r765}
\delta_{q,z} f(z)= \frac{{}f(q^{\frac{1}{2}}z)-f(q^{\frac{-1}{2}}z)}{(q^{\frac{1}{2}}-q^{\frac{-1}{2}})z}, \quad z\neq 0.
\end{equation}

 Jackson, see \cite{Jackson}, introduced the following three $q$-analogs of the Bessel function
 \begin{align}\label{s765}
\begin{split}
    J_{\alpha}^{(1)}(z;q)&=\frac{(q^{\alpha+1};q)_{\infty}}{(q;q)_{\infty}}\sum_{n=0}^{\infty}(-1)^n\frac{(\frac{z}{2})^{2n+\alpha}}{(q;q)_{n}(q^{\alpha+1};q)_{n}}
    \quad(|z|<2),\\
    J_{\alpha}^{(2)}(z;q)&=\frac{(q^{\alpha+1};q)_{\infty}}{(q;q)_{\infty}}\sum_{n=0}^{\infty}(-1)^n\frac{q^{n(\alpha+n)}(\frac{z}{2})^{2n+\alpha}}{(q;q)_{n}(q^{\alpha+1};q)_{n}}\quad(z\in \mathbb{C} ),\\
    J_{\alpha}^{(3)}(z;q)&=\frac{(q^{\alpha+1};q)_{\infty}}{(q;q)_{\infty}}\sum_{n=0}^{\infty}(-1)^n\frac{q^{\frac {n(n+1)}{2}}z^{2n+\alpha}}{(q;q)_{n}(q^{\alpha+1};q)_{n}}
    \quad(z\in \mathbb{C}).\end{split}
\end{align}
Set
\begin{equation}\label{g7907}
\begin{split}
\mathcal{J}_{\alpha}^{(k)} (z;q)=\left\{
  \begin{array}{ll}
    \dfrac{(q;q)_{\infty}}{(q^{\alpha+1};q)_{\infty}}(\frac {t}{2})^{-\alpha} J_{\alpha}^{(k)} (z;q) & \hbox{$( k=1,2)$;} \\
    \dfrac{(q;q)_{\infty}}{(q^{\alpha+1};q)_{\infty}}t^{-\alpha} J_{\alpha}^{(3)}(z;q)& \hbox{$(k=3)$.}
  \end{array}
\right.\end{split}
\end{equation}
  In some literature,  $\mathcal{J}_{\alpha}^{(k)} (z;q)\, (k=1,2,3)$  are called the modified Jackson $q$-Bessel functions, see \cite{modified}.

\vskip0.1cm

 Jackson, \cite{Jackson2}, introduced two $q$-analogs of the exponential function  by
 \begin{equation*}\label{h76}
\begin{split}
& E_q(z)= (-z(1-q);q)_{\infty}=\sum_{n=0}^{\infty}\frac{q^{\frac{n(n-1)}{2}}}{(q;q)_{n}}(z(1-q))^{n}, \,\, z\in \mathbb{C}, \\&
  e_q(z)=\frac{1}{(z(1-q);q)_{\infty}}=\sum_{n=0}^{\infty}\frac{(z(1-q))^{n}}{(q;q)_{n}}, \,\, |z|< \frac{1}{1-q}.
\end{split}
  \end{equation*}
 Exton, \cite{exton},  introduced a $q$-exponential function by
\begin{equation*}
\exp_{q}(z)=\sum_{n=0}^{\infty}\frac {q^{\frac{n(n-1)}{4}}}{(q;q)_{n}}(z(1-q))^{n},\quad z\in\mathbb{C}.
 \end{equation*}

The three $q$-pair  $\{\sin_{q}z,\,\, \cos z\},$   $\{Sin_{q}z,\,\, Cos_{q}z\}$, and  $\{S_{q}z,\,\, C_{q}z\}$ are defined  through the identities
\begin{align}\label{h765}
\begin{split}
  e_{q}(iz)&=\cos_{q}z+i\sin_{q}z,\\
E_{q}(iz)&=Cos_{q}z+iSin_{q}z,\\
\exp_{q}(iz)&=C_{q}z+iS_{q}z,\\\end{split}
\end{align}
respectively, see \cite{ Annaby,Gasper, Cardoss}.

\vskip 0.1 cm
  In this paper, we use the notations $\zeta_{1,q},\,\,\eta_{1,q}, \,\,\lambda_{1,q}$ and $\mu_{1,q}$  to denote  the smallest positive zero of the functions $Sin_{q}z,\,\,Cos_{q}z,\,\,S_{q}z$ and $C_{q}z$, respectively.

\vskip 0.1 cm

In \cite{sahar}, we defined the generalized $q$-Bernoulli polynomials $B^{(k)}_{n,\alpha}(z;q)\,\,( k=1,2,3)$ by the generating functions
 \begin{equation}\label{q14}
  \quad \frac{e_{q}(zt)e_{q}(\frac{-t}{2})}{g^{(1)}_{\alpha}(it;q)}=\sum_{n=0}^{\infty}B^{(1)}_{n,\alpha}(z;q)\frac{t^{n}}{[n]_{q}!},\quad|t|<
 \frac{ j^{(1)}_{1,\alpha}}{1-q},
\end{equation}
\begin{equation}\label{q13}
  \quad \frac{E_{q}(zt)E_{q}(\frac{-t}{2})}{g^{(2)}_{\alpha}(it;q)}=\sum_{n=0}^{\infty}B^{(2)}_{n,\alpha}(z;q)\frac{t^{n}}{[n]_{q}!},\quad|t|<\frac{
  j^{(2)}_{1,\alpha}}{1-q},
\end{equation}
\begin{equation}\label{q197}
   \frac{\exp_{q}(zt)\exp_{q}(\frac{-t}{2})}{g^{(3)}_{\alpha}(it;q)}=\sum_{n=0}^{\infty}B^{(3)}_{n,\alpha}(z;q)\frac{t^{n}}{[n]_{q}!},
   \quad|t|<\frac{2q^{\frac{1}{4}}j^{(3)}_{1,\alpha}}{1-q},
\end{equation}
where $j_{1,\alpha}^{(k)}$ is the first positive zero of $ J^{(k)}_{\alpha}(\cdot;q^{2}) \,\,( k=1,2,3)$, and  $g^{(k)}_{\alpha}(t;q) \,\,( k=1,2,3 )$  are the functions defined by
 \begin{equation*}
  g^{(k)}_{\alpha}(t;q)=\left\{
                          \begin{array}{ll}
                           \mathcal{J}_{\alpha}^{(k)} (t(1-q);q^{2})  & \hbox{$(k=1,2)$;} \\
                            \mathcal{J}_{\alpha}^{(3)}(\frac {t}{2}(1-q)q^{\frac{-1}{4}};q^{2}) & \hbox{$(k=3)$.}
                          \end{array}
                        \right.
\end{equation*}
We also defined the generalized $q$-Bernoulli numbers ${\beta}^{(k)}_{n,\alpha}(q)$ by
 $$\beta _{n,\alpha}^{(k)}:=B^{(k)}_{n,\alpha}(0;q)\,\,\,(k=1,2,3).$$
In \cite[Lemma 2.2]{sahar}, we proved
\begin{align}\label{m654}
 \frac{e_{q}(\frac{-t}{2})}{g^{(1)}_{\alpha}(it;q)}=\frac{E_{q}(\frac{-t}{2})}{g^{(2)}_{\alpha}(it;q)}.
\end{align}
 Consequently,
$$\beta_{n,\alpha}:={\beta}^{(1)}_{n,\alpha}(q)={\beta}^{(2)}_{n,\alpha}(q).$$

Hence
 \begin{equation}\label{q66}
 \frac{e_{q}(\frac{-t}{2})}{g^{(1)}_{\alpha}(it;q)}=\frac{E_{q}(\frac{-t}{2})}{g^{(2)}_{\alpha}(it;q)} =\sum_{n=0}^{\infty}\beta_{n,\alpha}(q)\frac{t^{n}}{[n]_{q}!},
\end{equation}
\begin{equation}\label{r66}
 \frac{\exp_{q}(\frac{-t}{2})}{g^{(3)}_{\alpha}(it;q)} =\sum_{n=0}^{\infty}\beta^{(3)}_{n,\alpha}(q)\frac{t^{n}}{[n]_{q}!}.
\end{equation}
\vskip0.1cm
The cases $\alpha=\pm\frac{1}{2}$ of Equations (\ref{q14})-(\ref{q197}) give the following $q$-Bernoulli and $q$-Euler polynomials introduced in \cite{IZ,marim},
 \begin{align}\label{q8}
 \begin{split}
  \frac{te_{q}(zt)}{e_{q}(\frac{t}{2})E_{q}(\frac{t}{2})-1}&=\sum_{n=0}^{\infty}b_{n}(z;q)\frac{t^{n}}{[n]_{q}!},\\
  \frac{t E_{q}(zt)}{e_{q}(\frac{t}{2})E_{q}(\frac{t}{2})-1}&=\sum_{n=0}^{\infty}B_{n}(z;q)\frac{t^{n}}{[n]_{q}!},\end{split}
\end{align}

\begin{align}\label{q13248}
 \begin{split}
  \frac{2e_{q}(zt)}{E_{q}(\frac{t}{2})e_{q}(\frac{t}{2})+1}&=\sum_{n=0}^{\infty}e_{n}(z;q)\frac{t^{n}}{[n]_{q}!},\\
  \frac{2 E_{q}(zt)}{E_{q}(\frac{t}{2})e_{q}(\frac{t}{2})+1}&=\sum_{n=0}^{\infty}E_{n}(z;q)\frac{t^{n}}{[n]_{q}!},\end{split}
\end{align}
and
\begin{equation}\label{q4039}
\begin{split}
  &\frac{t \exp_{q}(zt)\exp_{q}(\frac{-t}{2})}{\exp_{q}(\frac{t}{2})-\exp_{q}(\frac{-t}{2})}=\sum_{n=0}^{\infty}\tilde{B}_{n}(z;q)\frac{t^{n}}{[n]_{q}!},\\&
   \frac{2 \exp_{q}(zt)\exp_{q}(\frac{-t}{2})}{\exp_{q}(\frac{t}{2})+\exp_{q}(\frac{-t}{2})}=\sum_{n=0}^{\infty}\tilde{E}_{n}(z;q)\frac{t^{n}}{[n]_{q}!}.
   \end{split}
\end{equation}
\vskip0.1cm
The following identity
\begin{equation}\label{d8}
 g^{(1)}_{\alpha}(it;q)E_{q}(\frac{t}{2})= g^{(2)}_{\alpha}(it;q)e_{q}(\frac{t}{2})
={}_{2}\phi_{1}\,(q^{\alpha+\frac{1}{2}},-q^{\alpha+\frac{1}{2}};q^{2\alpha+1};q,\frac{(1-q)t}{2}),
\end{equation}
which we proved in \cite[Lemma 2.16]{sahar} will be needed in the the proof of some results of this paper.
\vskip0.2cm
The following Lemma from  \cite{marim} gives the reciprocal of $ \exp_{q}(z)$  in a certain domain.
\begin{lemma}\label{ii}
Let   $\displaystyle \Omega= \{ z\in \mathbb{C}: \, |1-exp_q(-z)|< 1\}$. Then
\begin{equation*}
\dfrac{1}{exp_q(z)}= \sum_{n=0}^{\infty} c_n\,z^{n}\quad (z\in \Omega),
\end{equation*}
where
\begin{equation}\label{a643}
 c_n=\, \sum_{k=1}^{n} (-1)^{k}\sum_{s_1+s_2+\ldots+s_k=n\atop  s_i>0\,(i=1,\ldots,k) } \dfrac{ q^{ \sum_{i=1}^k s_i(s_i-1)/4 }}{[s_1]_q! [s_2]_q!\ldots [s_{k}]_q!}.
 \end{equation}
\end{lemma}
\section{A determinant approach for the generalized $q$-Bernoulli polynomials}

This section introduces a determinant representation for the generalized $q$-Bernoulli polynomials defined in (\ref{q14})-(\ref{q197}). We begin by recalling some essential definitions and results  used in the sequel.
\vskip 0.1 cm
 Appell sequences are defined  for the first time in (1880) by  Appell \cite{Appell}, as
 any polynomial sequence $(p_{n}(\cdot))_{n}$ satisfying the identity
\begin{equation*}
  \frac{d}{dz}p_{n}(z)=np_{n-1}(z),\,\,p_{0}(z)\neq0.
\end{equation*}
 AL-Salam,  \cite{Alsalam1}, defined a $q$-analog of Appell sequence to be
 any set of polynomials $(p_{n}(\cdot))_{n}$  satisfying
\begin{equation}\label{d3}
  D_{q}p_{n}(z)=[n]_{q}p_{n-1}(z)\quad (n\in\mathbb{N}_{0}),
\end{equation}
where $[n]_{q}:=\dfrac{1-q^{n}}{1-q}$ is the $q$-natural number.  Keleshteri and  Mahmoudov, see \cite{two},  proved that the generating function of the  polynomials $(p_{n}(\cdot))_{n}$ defined in (\ref{d3}) is given by

\begin{equation}\label{t23}
  A_{q}(z;t)=A_{q}(t) e_{q}(zt)=\sum_{n=0}^{\infty}P_{n}(z)\frac{t^{n}}{[n]_{q}!},
\end{equation}
where
\begin{equation}\label{t24}
  A_{q}(t)=\sum_{n=0}^{\infty}A_{n,q}\frac{t^{n}}{[n]_{q}!},\,\,A_{q}(t)\neq0,
\end{equation}
is an analytic function at $t = 0$ and  $A_{n,q}= P_{n}(0)$.

\vskip 0.3 cm

In \cite{Sadj},  Njionou called the sequence of polynomials  that satisfy
 $$D_{q}P_{n}(z)=[n]_{q}P_{n-1}(qz),$$
  or equivalently
\begin{equation}\label{d63}
  D_{q^{-1}}P_{n}(z)=[n]_{q}P_{n-1}(z)\quad (n\in\mathbb{N}),
\end{equation}
 the $q$-Appell polynomials of type II. He also proved that the polynomials defined by (\ref{d63}) have the generating function

 \begin{equation}\label{t273}
  A_{q}(z;t)=A_{q}(t) E_{q}(zt)=\sum_{n=0}^{\infty}P_{n}(z)\frac{t^{n}}{[n]_{q}!},
\end{equation}
 where

\begin{equation}\label{t284}
  A_{q}(t)=\sum_{n=0}^{\infty}A_{n,q}\frac{t^{n}}{[n]_{q}!},\,\,A_{q}(t)\neq0,
\end{equation}

is an analytic function at $t = 0$ and $A_{n,q}= P_{n}(0)$.

\vskip 0.2 cm

 The following theorem is proved in \cite{two}.

\begin{theorem}\label{yy}
Suppose that $\{P_{n}(z)\}$  are the sequence of $q$-Appell polynomials with
generating function $A_{q}(z;t)$, defined in the relations (\ref{t23}) and (\ref{t24}). If
 $$\dfrac{1}{A_{q}(t)}=\sum_{k=0}^{\infty}\mu_{k,q}\frac{t^{k}}{k!},$$ then
 $P_{0}(z)=\dfrac{1}{\mu_{0,q}}$ and for $n\in\mathbb{N}$

 \begin{equation}\label{t12089}
\begin{split}
 P_{n}(z)&=\frac{(-1)^{n}}{(\mu_{0,q})^{n+1}}\\&\left|
     \begin{array}{ccccccc}
       1 & z & z^{2} & \ldots & \ldots & z^{n-1} & z^{n} \\
       \mu_{0,q} & \mu_{1,q} & \mu_{2,q} & \ldots& \ldots & \mu_{n-1,q} & \mu_{n,q} \\
       0 & \mu_{0,q} &  \left[ \begin{array}{c}
                                                      2 \\
                                                      1
                                                    \end{array}\right]_{q}\mu_{1,q}& \ldots & \ldots & \left[ \begin{array}{c}
                                                      n-1 \\
                                                      1
                                                    \end{array}\right]_{q}\mu_{n-2,q} & \left[ \begin{array}{c}
                                                      n \\
                                                      1
                                                    \end{array}\right]_{q}\mu_{n-1,q} \\
      0 & 0 & \mu_{0,q} & \ldots & \ldots & \left[ \begin{array}{c}
                                                      n-1 \\
                                                      2
                                                    \end{array}\right]_{q}\mu_{n-3,q} & \left[ \begin{array}{c}
                                                      n \\
                                                      2
                                                    \end{array}\right]_{q}\mu_{n-2,q} \\
       \vdots &  &  & \ddots &  & \vdots & \vdots \\
       \vdots &  &  &  & \ddots & \vdots & \vdots \\
       0 & \ldots & \ldots & \ldots & 0 & \mu_{0,q} &\left[ \begin{array}{c}
                                                      n \\
                                                      n-1
                                                    \end{array}\right]_{q}\mu_{1,q}
        \\
     \end{array}
   \right|.
   \end{split}
\end{equation}

\end{theorem}

Under the conventions that $\left[
                              \begin{array}{c}
                                0 \\
                                0 \\
                              \end{array}
                            \right]_{q}=1,\,\,\left[
                              \begin{array}{c}
                                n \\
                                k \\
                              \end{array}
                            \right]_{q}=0,\,\,\text {for}\,\,\,k>n\geq 0$,
the polynomials $P_{n}(z)$ defined in (\ref{t12089}) can be written in the form

 \begin{equation*}
\begin{split}
 P_{n}(z)&=\dfrac{(-1)^{n}}{(\mu_{0,q})^{n+1}}\,\text {det}\, A_{n}\quad (n\in\mathbb{N}),
   \end{split}
\end{equation*}
where $A_{n}=(a_{ij})^{n}_{i,j=0}$ is the $(n+1)\times (n+1)$ matrix defined for $j\in\{0,1,\cdots,n\}$ by
 \begin{equation*}
 a_{ij}=\left\{
  \begin{array}{ll}
   z^{j}, & \hbox{$i=0$;}\\
 \\
   \left[
     \begin{array}{c}
       j \\
       i-1 \\
     \end{array}
   \right]_{q}
 \mu_{j-i+1,q}, & \hbox{$1\leq i\leq n$.}
  \end{array}
\right.
 \end{equation*}

\begin{theorem}\label{rt}
The generalized $q$-Bernoulli polynomials $(B^{(1)}_{n,\alpha}(z;q))_{n}$ have the representation

\begin{equation}\label{t594}
\begin{split}
B^{(1)}_{n,\alpha}(z;q)=
    (-1)^{n}\,\text {det}\, A_{n}  \quad (n\in\mathbb{N}_{0}),
   \end{split}
\end{equation}

where $A_{n}=(a_{ij})^{n}_{i,j=0}$ is the $(n+1)\times(n+1)$ matrix defined for $j\in\{0,1,\cdots,n\}$ by
 \begin{equation}\label{d7}
 a_{ij}=\left\{
  \begin{array}{ll}
   z^{j}, & \hbox{$i=0$;}\\
 \\
   \left[
     \begin{array}{c}
       j \\
       i-1 \\
     \end{array}
   \right]_{q}
 \dfrac{1}{2^{j-i+1}}\dfrac{(q^{2\alpha+1};q^{2})_{j-i+1}}{(q^{2\alpha+1};q)_{j-i+1}}, & \hbox{$1\leq i\leq n$.}
  \end{array}
\right.
 \end{equation}
\end{theorem}
\vskip 0.1 cm
\begin{proof}
Since the polynomials $(B^{(1)}_{n,\alpha}(z;q))_{n}$ satisfy the $q$-difference equation (\ref{d3}),  then,  from (\ref{t23}), the generating function of $B^{(1)}_{n,\alpha}(z;q)$
can be written as
\begin{equation*}
  A_{q}(t)e_{q}(zt)=\sum_{n=0}^{\infty}B^{(1)}_{n,\alpha}(z;q)\frac{t^{n}}{[n]_{q}!},\quad|t|<\frac {j^{(1)}_{1, \alpha}}{1-q},
\end{equation*}
where
\begin{equation*}
  A_{q}(t)= \frac{e_{q}(\frac{-t}{2})}{g^{(1)}_{\alpha}(it;q)}.
\end{equation*}
Now set
\begin{equation*}
 B_{q}(t):=\frac{1}{A_{q}(t)}=\frac{g^{(1)}_{\alpha}(it;q)}{e_{q}(\frac{-t}{2})}=g^{(1)}_{\alpha}(it;q)E_{q}(\frac{t}{2}).
\end{equation*}
Thus from (\ref{d8}),
\begin{equation}\label{f17981}
\begin{split}
 B_{q}(t)= \sum_{n=0}^{\infty}\frac{(q^{2\alpha+1};q^{2})_{n}}{2^{n}(q^{2\alpha+1};q)_{n}}\frac{t^{n}}{[n]_{q}!}
=\sum_{n=0}^{\infty}\mu_{n,\alpha}(q)\frac{t^{n}}{[n]_{q}!}.
 \end{split}
\end{equation}
 Equating the coefficients of $t^{n}$ in (\ref{f17981}), we obtain

\begin{equation*}
\begin{split}
  \mu_{n,\alpha}(q)=\dfrac{1}{2^{n}}\dfrac{(q^{2\alpha+1};q^{2})_{n}}{(q^{2\alpha+1};q)_{n}}\quad\,\, (n\in\mathbb{N}_{0}).
  \end{split}
\end{equation*}

Therefore, from Theorem \ref{yy},  $ B^{(1)}_{0,\alpha}(z;q)=1$ and

\begin{equation*}
\begin{split}
B^{(1)}_{n,\alpha}(z;q) =(-1)^{n}\,\text {det}\, A_{n}\quad (n\in\mathbb{N}),
   \end{split}
\end{equation*}
where $A_{n}$ is the matrix defined in (\ref{d7}).
\end{proof}
\vskip 0.1 cm
\begin{corollary}
The generalized $q$-Bernoulli numbers  $(\beta_{n,\alpha}(q))_{n}$ defined in (\ref{q66}) are given by $\beta_{0,\alpha}(q)=1$ and

\begin{equation*}
\begin{split}
\beta_{n,\alpha}(q) =(-1)^{n}\,\text {det}\, A_{n}\quad (n\in\mathbb{N}),
   \end{split}
\end{equation*}
where $A_{n}=(a_{ij})^{n}_{i,j=1}$ is the $n\times n$ matrix defined by
 \begin{equation*}
 a_{ij}=\
   \left[
     \begin{array}{c}
       j \\
       i-1 \\
     \end{array}
   \right]_{q}
 \dfrac{1}{2^{j-i+1}}\dfrac{(q^{2\alpha+1};q^{2})_{j-i+1}}{(q^{2\alpha+1};q)_{j-i+1}}, \quad 1\leq i,j\leq n.
 \end{equation*}
\end{corollary}
\vskip 0.1 cm
\begin{proof}
The proof follows directly by setting $z=0$ in Theorem \ref{rt}.
\end{proof}
\begin{theorem}\label{js}
Let  $\{P_{n}(z)\}$  be the $q$-Appell polynomials defined in (\ref{t273}). If  $$\dfrac{1}{A_{q}(t)}=\displaystyle\sum_{k=0}^{\infty}\mu_{k,q}\frac{t^{k}}{k!},$$ then
\begin{equation*}
 P_{n}(z)=
    \dfrac{(-1)^{n}}{(\mu_{0,q})^{n+1}}\text {det}\, A_{n}\quad (n\in\mathbb{N}_{0}),
\end{equation*}

where $A_{n}=(a_{ij})^{n}_{i,j=0}$ is the $(n+1)\times (n+1)$ matrix defined for $j\in\{0,1,\cdots,n\}$ by
 \begin{equation*}
 a_{ij}=\left\{
  \begin{array}{ll}
   q^{\frac{j(j-1)}{2}}z^{j}, & \hbox{$i=0$;}\\
 \\
   \left[
     \begin{array}{c}
       j \\
       i-1 \\
     \end{array}
   \right]_{q}
 \mu_{j-i+1,q}, & \hbox{$1\leq i\leq n$.}
  \end{array}
\right.
 \end{equation*}

\end{theorem}

\vskip 0.1 cm

\begin{proof}
 The proof is similar to the proof of Theorem \ref{yy} in \cite{two} and is omitted.
\end{proof}
\begin{theorem}\label{ht}
The generalized $q$-Bernoulli polynomials $(B^{(2)}_{n,\alpha}(z;q))_{n}$ are represented by

\begin{equation}\label{h13}
\begin{split}
B^{(2)}_{n,\alpha}(z;q)=
    (-1)^{n}\,\text {det}\, A_{n}\quad (n\in\mathbb{N}_{0}),
   \end{split}
\end{equation}

where $A_{n}=(a_{ij})^{n}_{i,j=0}$ is the $(n+1)\times (n+1)$ matrix defined for $j\in\{0,1,\cdots,n\}$ by
 \begin{equation}\label{d9}
 a_{ij}=\left\{
  \begin{array}{ll}
   q^{\frac {j(j-1)}{2}}z^{j}, & \hbox{$i=0$;}\\
 \\
   \left[
     \begin{array}{c}
       j \\
       i-1 \\
     \end{array}
   \right]_{q}
 \dfrac{1}{2^{j-i+1}}\dfrac{(q^{2\alpha+1};q^{2})_{j-i+1}}{(q^{2\alpha+1};q)_{j-i+1}}, & \hbox{$1\leq i\leq n$.}
  \end{array}
\right.
 \end{equation}

\end{theorem}

\vskip 0.1 cm

\begin{proof}
Since the generalized $q$-Bernoulli polynomials $(B^{(2)}_{n,\alpha}(z;q))_{n}$  satisfy the $q$-difference equation (\ref{d63}),
then,  from (\ref{t273}), the generating function of $B^{(2)}_{n,\alpha}(z;q)$ can be written as
\begin{equation*}
  A_{q}(t)E_{q}(zt)=\sum_{n=0}^{\infty}B^{(2)}_{n,\alpha}(z;q)\frac{t^{n}}{[n]_{q}!},\quad |t|<\frac {j_{1,\alpha}^{(2)}}{1-q},
\end{equation*}
where
\begin{equation*}
  A_{q}(t)=  \frac{E_{q}(\frac{-t}{2})}{g^{(2)}_{\alpha}(it;q)}.
\end{equation*}
Set
\begin{equation*}
  B_{q}(t):=\frac{1}{A_{q}(t)}=\frac{g^{(2)}_{\alpha}(it;q)}{E_{q}(\frac{-t}{2})}=g^{(2)}_{\alpha}(it;q)e_{q}(\frac{t}{2}).
\end{equation*}
Hence from (\ref{d8}),
\begin{equation}\label{g41}
  B_{q}(t)= \sum_{n=0}^{\infty}\frac{(q^{2\alpha+1};q^{2})_{n}}{2^{n}(q^{2\alpha+1};q)_{n}}\frac{t^{n}}{[n]_{q}!}=\sum_{n=0}^{\infty}\mu_{n,\alpha}(q)\frac{t^{n}}{[n]_{q}!}.
\end{equation}
So, equating the coefficients of $t^{n}$ in (\ref{g41}) gives
\begin{equation*}\label{t33}
\begin{split}
  \mu_{n,\alpha}(q)=\dfrac{1}{2^{n}}\dfrac{(q^{2\alpha+1};q^{2})_{n}}{(q^{2\alpha+1};q)_{n}} \quad (n\in\mathbb{N}_{0}).\end{split}
\end{equation*}
Consequently, by applying Theorem \ref{js}, we obtain the identities in (\ref{h13}), (\ref{d9}) and complete the proof.
\end{proof}

\vskip 0.1 cm
\begin{remark}
 We can  obtain the result of Theorem \ref{ht} by replacing $q$ by $q^{-1}$ in Theorem \ref{rt} and using
\begin{equation}\label{e3}
   B^{(2)}_{n,\alpha}(z;q)=q^{\frac{n(n-1)}{2}} B^{(1)}_{n,\alpha} (z;\frac{1}{q}),
\end{equation}
see \cite[Proposition 2.12]{sahar}. Indeed, replacing $q $ by $q^{-1}$  in (\ref{t594}) yields

\begin{equation}\label{e7}
\begin{split}
B^{(1)}_{n,\alpha}(z;\frac{1}{q}) =(-1)^{n}\,det\, A_{n}\quad (n\in\mathbb{N}_{0}),
   \end{split}
\end{equation}
where $A_{n}=(a_{ij})^{n}_{i,j=0}$ is the $(n+1)\times(n+1)$ matrix defined for $j\in\{0,1,\cdots,n\}$ by
 \begin{equation*}
 a_{ij}=\left\{
  \begin{array}{ll}
   z^{j}, & \hbox{$i=0$;}\\
 \\
   \left[
     \begin{array}{c}
       j \\
       i-1 \\
     \end{array}
   \right]_{q^{-1}}
 \dfrac{1}{2^{j-i+1}}\dfrac{(q^{-2\alpha-1};q^{-2})_{j-i+1}}{(q^{-2\alpha-1};q^{-1})_{j-i+1}}, & \hbox{$1\leq i\leq n$.}
  \end{array}
\right.
 \end{equation*}
Using the identities, cf. \cite[Eq. I.3, I.47]{Gasper},
\begin{equation*}
  (a;q^{-1})_{n}=(a^{-1};q)_{n}\,(-a)^{n}q^{-\frac{n(n-1)}{2}}\quad\text{ and }\quad\left[
                                                                                        \begin{array}{c}
                                                                                          \alpha \\
                                                                                          k \\
                                                                                        \end{array}
                                                                                      \right]_{q^{-1}}=q^{k^{2}-\alpha k}\left[
                                                                                        \begin{array}{c}
                                                                                          \alpha \\
                                                                                          k \\
                                                                                        \end{array}
                                                                                      \right]_{q},
\end{equation*}
we obtain for $j\in\{0,1,\cdots,n\}$
 \begin{equation*}
 \begin{split}
 a_{ij}=\left\{
  \begin{array}{ll}
   z^{j}, & \hbox{$i=0$;}\\
 \\
   \left[
     \begin{array}{c}
       j \\
       i-1 \\
     \end{array}
   \right]_{q}
 q^{\frac{(i-1)(i-2)-j(j-1)}{2}}\dfrac{(q^{2\alpha+1};q^{2})_{j-i+1}}{2^{j-i+1}(q^{2\alpha+1};q)_{j-i+1}}, & \hbox{$1\leq i\leq n$.}
  \end{array}
\right.
\end{split}
 \end{equation*}
 Taking $q^{\frac{-j(j-1)}{2}}$  from the $j^{th}$ column of $A_{n}$ as a common factor, $j\in\{0,1,\cdots,n\}$ gives
 \begin{equation*}
 \text {det}\, A_{n}=\left(\prod_{j=0}^{n}q^{\frac{-j(j-1)}{2}}\right)\,\text{det}\, (B_{n}),
 \end{equation*}
where $B_{n}=(b_{ij})_{i,j=0}^{n},$
 \begin{equation*}
 b_{ij}=\left\{
  \begin{array}{ll}
   q^{\frac{j(j-1)}{2}}z^{j}, & \hbox{$i=0$;}\\
 \\
   \left[
     \begin{array}{c}
       j \\
       i-1 \\
     \end{array}
   \right]_{q}
 q^{\frac{(i-1)(i-2)}{2}}\dfrac{(q^{2\alpha+1};q^{2})_{j-i+1}}{2^{j-i+1}(q^{2\alpha+1};q)_{j-i+1}}, & \hbox{$1\leq i\leq n$.}
  \end{array}
\right.
 \end{equation*}
 Now, taking $q^{\frac{(i-1)(i-2)}{2}}$ as a common factor from the $i^{th}$ row of $B_{n}, \,i\in\{2,3,\cdots,n\}$  yields
   \begin{equation*}
    \begin{split}
 \text {det}\, A_{n}=\left(\prod_{j=0}^{n}q^{\frac{-j(j-1)}{2}}\right)\left( \prod_{i=2}^{n}q^{\frac{(i-1)(i-2)}{2}}\right)\,\text{det}\, (c_{n})=q^{\frac{-n(n-1)}{2}}\,\text{det}\, (c_{n}),
  \end{split}
 \end{equation*}
where $c_{n}=(c_{ij})_{i,j=0}^{n},$

 \begin{equation*}
 c_{ij}=\left\{
  \begin{array}{ll}
   q^{\frac{j(j-1)}{2}}z^{j}, & \hbox{$i=0$;}\\
 \\
   \left[
     \begin{array}{c}
       j \\
       i-1 \\
     \end{array}
   \right]_{q}
\dfrac{1}{2^{j-i+1}} \dfrac{(q^{2\alpha+1};q^{2})_{j-i+1}}{(q^{2\alpha+1};q)_{j-i+1}}, & \hbox{$1\leq i\leq n$.}
  \end{array}
\right.
 \end{equation*}

 Hence
 \begin{equation}\label{e4}
\begin{split}
       B^{(2)}_{n,\alpha}(z;q)=q^{\frac{n(n-1)}{2}}B^{(1)}_{n,\alpha}(z;\frac{1}{q}) =(-1)^{n}\,\text {det}\, (c_{n})\quad (n\in\mathbb{N}),
   \end{split}
\end{equation}
and this yields (\ref{h13}).
\end{remark}
\vskip0.4cm

In the following, we introduce a determinant representation for the generalized $q$-Bernoulli polynomials $B^{(3)}_{n,\alpha}(z;q)\,\,(n\in\mathbb{N}_{0})$  by the same technique we used for $B^{(k)}_{n,\alpha}(z;q)\,\,(k=1,2)$.  The  polynomials $B^{(3)}_{n,\alpha}(z;q)$ satisfy the $q$-difference equation
\begin{equation}\label{d87}
 \delta_{q,z}\,B^{(3)}_{n,\alpha}(z;q)=[n]_{q}B^{(3)}_{n-1,\alpha}(z;q)\quad (n\in\mathbb{N}),
\end{equation}
where $\delta_{q,z}$ is the operator defined in (\ref{r765}).  Hence the generating function of $B^{(3)}_{n,\alpha}(z;q)$ can be written as

\begin{equation}\label{t1543}
 A_{q}(t)\exp_{q}(zt)=\sum_{n=0}^{\infty}B^{(3)}_{n,\alpha}(z;q)\frac{t^{n}}{[n]_{q}!},\quad|t|<\frac{2q^{\frac{1}{4}}j^{(3)}_{1,\alpha}}{1-q},
\end{equation}
where
\begin{equation}\label{t9088}
  A_{q}(t)=\frac{\exp_{q}(\frac{-t}{2})}{g^{(3)}_{\alpha}(it;q)},\,\,\,A_{q}(t)\neq0,
\end{equation}
is an analytic function at $t = 0$, and $j^{(3)}_{1,\alpha}$ is the smallest positive zero of $J_{\alpha}^{(3)}(\cdot;q^{2})$.

\vskip 0.3 cm
\begin{theorem}\label{nn}
Let $\{P_{n}(z)\}$  be the $q$-Appell polynomials defined in (\ref{t1543}). If  $$\dfrac{1}{A_{q}(t)}=\displaystyle\sum_{k=0}^{\infty}\mu_{k,q}\frac{t^{k}}{k!},$$ then

 \begin{equation*}
\begin{split}
P_{n}(z)=
   \dfrac{(-1)^{n}}{(\mu_{0,q})^{n+1}}\,\text {det}\, A_{n} \quad (n\in\mathbb{N}_{0}),
   \end{split}
\end{equation*}
where $A_{n}=(a_{ij})^{n}_{i,j=0}$ is the $(n+1)\times( n+1)$ matrix defined for $j\in\{0,1,\cdots,n\}$ by
 \begin{equation*}
 a_{ij}=\left\{
  \begin{array}{ll}
   q^{\frac{j(j-1)}{4}}z^{j}, & \hbox{$i=0$;}\\
 \\
   \left[
     \begin{array}{c}
       j \\
       i-1 \\
     \end{array}
   \right]_{q}
 \mu_{j-i+1,q}, & \hbox{$1\leq i\leq n$.}
  \end{array}
\right.
 \end{equation*}
\end{theorem}

\vskip 0.1 cm
\begin{proof}
 The proof is similar to the proof of Theorem \ref{yy} in \cite{two} and is omitted.

\end{proof}

\begin{theorem}\label{st}
The generalized $q$-Bernoulli polynomial $(B^{(3)}_{n,\alpha}(z;q))_{n}$  are represented by
\begin{equation*}
\begin{split}
B^{(3)}_{n,\alpha}(z;q)=
    (-1)^{n}\,\text {det}\, A_{n} \quad (n\in\mathbb{N}_{0}),
   \end{split}
\end{equation*}

where $A_{n}=(a_{ij})^{n}_{i,j=0}$ is the $(n+1)\times (n+1)$ matrix defined for $j\in\{0,1,\cdots,n\}$ by
 \begin{equation*}
 a_{ij}=\left\{
  \begin{array}{ll}
   q^{\frac {j(j-1)}{4}}z^{j}, & \hbox{$i=0$;}\\
 \\
   \left[
     \begin{array}{c}
       j \\
       i-1 \\
     \end{array}
   \right]_{q}
 \mu_{j-i+1,\alpha}(q), & \hbox{$1\leq i\leq n$,}
  \end{array}
\right.
 \end{equation*}
where
\begin{equation}\label{h87}
\begin{split}
  \mu_{m,\alpha}(q)=\frac{(-1)^{m}[m]_{q}!}{2^{m}}\sum_{k=0}^{[\frac{m}{2}]}\frac{q^{k(k+\frac{1}{2})}(1-q)^{2k}c_{m-2k}
}{(q^{2},q^{2\alpha+2};q^{2})_{k}}\quad (m\in\mathbb{N}_{0}),\end{split}
\end{equation}
and the coefficients $(c_{k})_{k}$ are  defined in (\ref{a643}).

\end{theorem}
\vskip0.1cm
\begin{proof}
The generalized $q$-Bernoulli polynomials  $B^{(3)}_{n,\alpha}(z;q)$ are defined by (\ref{t1543}) and (\ref{t9088}). Now, set
\begin{equation*}
 B_{q}(t):=\frac{1}{A_{q}(t)}=\frac{g^{(3)}_{\alpha}(it;q)}{\exp_{q}(\frac{-t}{2})} = \displaystyle\sum_{m=0}^{\infty}\mu_{m,\alpha}(q)\frac{t^{m}}{[m]_{q}!}.
\end{equation*}
Then $B_{q}(t)A_{q}(t)=1$,  and by Lemma \ref{ii}, we get
\begin{equation*}
\begin{split}
  B_{q}(t)= g^{(3)}_{\alpha}(it;q)\frac{1}{\exp_{q}(\frac{-t}{2})}
  &=\left(\sum_{n=0}^{\infty}\frac{q^{n(n+\frac{1}{2})}(1-q)^{2n}t^{2n}}{2^{2n}(q^{2},q^{2\alpha+2};q^{2})_{n}}\right)
\left(\sum_{n=0}^{\infty}c_{n}\frac{(-t)^{n}}{2^{n}}\right)\\&=\sum_{m=0}^{\infty}\frac{(-1)^{m}t^{m}}{2^{m}}\sum_{k=0}^{[\frac{m}{2}]}\frac{q^{k(k+\frac{1}{2})}(1-q)^{2k}
c_{m-2k}
}{(q^{2},q^{2\alpha+2};q^{2})_{k}}.\end{split}
\end{equation*}
Consequently,
\begin{equation}\label{g40}
\begin{split}
  B_{q}(t)&=\sum_{m=0}^{\infty}\frac{(-1)^{m}t^{m}}{2^{m}[m]_{q}!}\sum_{k=0}^{[\frac{m}{2}]}\frac{[m]_{q}!q^{k(k+\frac{1}{2})}(1-q)^{2k}c_{m-2k}
}{(q^{2},q^{2\alpha+2};q^{2})_{k}}=\sum_{m=0}^{\infty}\mu_{m,\alpha}(q)\frac{t^{m}}{[m]_{q}!}.\end{split}
\end{equation}
Equating the coefficients of $t^{m}$ in the right hand  side of (\ref{g40}) gives (\ref{h87}). Hence by applying Theorem \ref{nn}, we obtain the desired result.

\end{proof}
\begin{corollary}
 The generalized $q$-Bernoulli numbers $(\beta^{(3)}_{n,\alpha}(q))_{n}$ defined in (\ref{r66}) are given by
\begin{equation*}
\begin{split}
\beta^{(3)}_{n,\alpha}(q) =(-1)^{n}\,\text {det}\, A_{n}\quad (n\in\mathbb{N}),
   \end{split}
\end{equation*}
where $A_{n}=(a_{ij})^{n}_{i,j=1}$ is the $n\times n$ matrix defined by
 \begin{equation*}
 a_{ij}=
   \left[
     \begin{array}{c}
       j \\
       i-1 \\
     \end{array}
   \right]_{q}
 \mu_{j-i+1,\alpha}(q),\quad 1\leq i,j\leq n,
 \end{equation*}
and $\mu_{m,\alpha}(q)\,\,(m\in\mathbb{N}_{0})$ is defined in (\ref{h87}).
\end{corollary}

\vskip 0.1 cm

\begin{proof}
The proof follows by setting $z=0$ in Theorem \ref{st}.
\end{proof}

\section{Asymptotic behavior of the generalized $q$-Bernoulli polynomials for sufficiently large $n$}
 In this section,  we use the  Darboux Method to determine the asymptotic behavior of the generalized $q$-Bernoulli polynomials for sufficiently large $n$.  The Darboux method is described as follows,  see \cite{olver}.
Let $r$ be the distance from the origin of the nearest singularity of $ f(t)=\displaystyle\sum_{n=0}^{\infty}a_{n}\, t^{n}$ and suppose that we can find a comparison function g(t) having the following properties

\vskip 0.4 cm

(i) $g(t)$ is analytic in $0<|t|<r$,

\vskip 0.2 cm
(ii) $f(t)-g(t)$  is continuous in $0<|t|\leq r$,

\vskip 0.2 cm
(iii) the coefficients $(b_{n})_{n}$ in the Laurent expansion

\begin{equation*}
  g(t)=\sum_{n=-\infty}^{\infty}b_{n}\, t^{n},\quad 0<|t|<r,
\end{equation*}
have known asymptotic behavior. Then
\begin{equation*}
  a_{n}=b_{n}+o(r^{-n})\quad (n\rightarrow\infty).
\end{equation*}

\begin{theorem}\label{rr}
The large $n$ asymptotic of the generalized $q$-Bernoulli polynomials $(B^{(2)}_{n,\alpha}(z;q))_{n}$ are
\begin{equation*}
\begin{split}
B^{(2)}_{2n,\alpha}(z;q)&=\dfrac{2}{1-q}\dfrac{(-1)^{n+1}[2n]_{q}!
\left(Cos_{q}(2\zeta_{\alpha,1}z)Cos_{q}(\zeta_{\alpha,1})+Sin_{q}(2\zeta_{\alpha,1}z)Sin_{q}(\zeta_{\alpha,1})
\right)}{(2\zeta_{\alpha,1})^{2n+1}\,\frac{d}{dz}\mathcal{J}^{(2)}_{\alpha}(z;q^{2})|_{z=j^{(2)}_{1,\alpha}}
}\\&+o(r^{-2n}),
\end{split}
\end{equation*}

\begin{equation*}
\begin{split}
B^{(2)}_{2n+1,\alpha}(z;q)&=\dfrac{2}{1-q}\dfrac{(-1)^{n+1}[2n+1]_{q}!\left(Sin_{q}(2\zeta_{\alpha,1}z)Cos_{q}(\zeta_{\alpha,1})
-Cos_{q}(2\zeta_{\alpha,1}z)Sin_{q}(\zeta_{\alpha,1})\right)}{(2\zeta_{\alpha,1})^{2n+2}\,\frac{d}{dz}\mathcal{J}^{(2)}_{\alpha}(z;q^{2})|_{z=j^{(2)}_{1,\alpha}}
}\\&+o(r^{-2n}),
\end{split}
\end{equation*}
as  $n\,\rightarrow\,\infty$, where the functions $Sin_{q}(\cdot)$ and $Cos_{q}(\cdot)$  are defined in (\ref{h765}),  $0<r\leq 2\zeta_{\alpha,1},$ and $\zeta_{\alpha,1}:=\dfrac{j^{(2)}_{1,\alpha}}{2(1-q)}$.

\end{theorem}

\vskip 0.1 cm

\begin{proof}
We prove the theorem by the Darboux method. Let $f(z,t)$ denote the left hand side of (\ref{q13}).  Then
\begin{equation*}
  f(z,t):=\frac{E_{q}(zt)E_{q}(\frac{-t}{2})}{g^{(2)}_{\alpha}(it;q)},\quad
  |t|<\frac{j^{(2)}_{1,\alpha}}{1-q}.
\end{equation*}
 Set
\begin{equation*}
\begin{split}
  g(z,t)&:=Res(f,\frac{ij^{(2)}_{1,\alpha}}{1-q})\frac{1}{(t-\frac {ij^{(2)}_{1,\alpha}}{1-q})}+Res(f,\frac{-ij^{(2)}_{1,\alpha}}{1-q})\frac{1}{(t+\frac {ij^{(2)}_{1,\alpha}}{1-q})},
\end{split}
\end{equation*}
\begin{equation}\label{r56}
  f(z,t)= \displaystyle\sum_{n=0}^{\infty}a_{n}\, t^{n}, \quad \text {and } \quad   g(z,t)=\displaystyle\sum_{n=0}^{\infty}b_{n}\, t^{n}.
\end{equation}
 Then
\begin{equation*}\label{y4}
\begin{split}
  g(z,t)=Res(f,2i\zeta_{\alpha,1})\frac{1}{(t-2i\zeta_{\alpha,1})} +Res(f,-2i\zeta_{\alpha,1})\frac{1}{(t+2i\zeta_{\alpha,1})}.
\end{split}
\end{equation*}
Clearly $f(z,t)-g(z,t)$ is continuous within and on the the circle $|t|=2\zeta_{\alpha,1}$.  Now, we calculate the coefficient of $t^{n}$
in the expansion of $g(z,t)$. Indeed,

\begin{equation}\label{y7}
\begin{split}
  g(z,t)&=Res(f,2i\zeta_{\alpha,1})\frac{1}{(t-2i\zeta_{\alpha,1})} +Res(f,-2i\zeta_{\alpha,1})\frac{1}{(t+2i\zeta_{\alpha,1})}
\\&=Res(f,2i\zeta_{\alpha,1})\frac{-1}{2i\zeta_{\alpha,1}}\sum_{n=0}^{\infty}(\frac{t}{2i\zeta_{\alpha,1}})^{n}+Res(f,-2i\zeta_{\alpha,1})\frac{1}{2i\zeta_{\alpha,1}}
\sum_{n=0}^{\infty}(\frac{-t}{2i\zeta_{\alpha,1}})^{n}.
\end{split}
\end{equation}
Since

\begin{equation*}
\begin{split}
&Res(f,2i\zeta_{\alpha,1} )=\frac{E_{q}(2i\zeta_{\alpha,1}z)E_{q}(-i\zeta_{\alpha,1})}{\frac{d}{dz}\mathcal{J}^{(2)}_{\alpha}(iz(1-q);q^{2})|_{z=2i\zeta_{\alpha,1}}},\\&
Res(f,-2i\zeta_{\alpha,1} )=\frac{E_{q}(-2i\zeta_{\alpha,1}z)E_{q}(i\zeta_{\alpha,1})}
{\frac{d}{dz}\mathcal{J}^{(2)}_{\alpha}(iz(1-q);q^{2})|_{z=-2i\zeta_{\alpha,1}}},
 \end{split}
\end{equation*}
then Equation (\ref{y7}) can be written as

\begin{equation*}
\begin{split}
  g(z,t)=&
-\frac{E_{q}(2i\zeta_{\alpha,1}z)E_{q}(-i\zeta_{\alpha,1})}{\frac{d}{dz}\mathcal{J}^{(2)}_{\alpha}(iz(1-q);q^{2})|_{z=2i\zeta_{\alpha,1}}}
\sum_{n=0}^{\infty}\frac{i^{-n-1}}{(2\zeta_{\alpha,1})^{n+1}}t^{n}
\\&-\frac{E_{q}(-2i\zeta_{\alpha,1}z)E_{q}(i\zeta_{\alpha,1})}
{\frac{d}{dz}\mathcal{J}^{(2)}_{\alpha}(iz(1-q);q^{2})|_{z=-2i\zeta_{\alpha,1}}}\sum_{n=0}^{\infty}\frac{(-i)^{-n-1}}{(2\zeta_{\alpha,1})^{n+1}}t^{n}.
\end{split}
\end{equation*}

Since
\begin{equation*}
\begin{split}
\frac{d}{dz}\mathcal{J}^{(2)}_{\alpha}(iz(1-q);q^{2})|_{z=2i\zeta_{\alpha,1}}&=i(1-q)\frac{d}{dz}\mathcal{J}^{(2)}_{\alpha}(iz(1-q);q^{2})|_{z=2i\zeta_{\alpha,1}}
\\&=-i(1-q)\frac{d}{dz}\mathcal{J}^{(2)}_{\alpha}(z;q^{2})|_{z=j^{(2)}_{1,\alpha}},
\end{split}
\end{equation*}
then

\begin{equation}\label{y787}
\begin{split}
  g(z,t)=&
-\frac{E_{q}(2i\zeta_{\alpha,1}z)E_{q}(-i\zeta_{\alpha,1})}{(1-q)\frac{d}{dz}\mathcal{J}^{(2)}_{\alpha}(z;q^{2})|_{z=j^{(2)}_{1,\alpha}}}
\sum_{n=0}^{\infty}\frac{i^{-n}}{(2\zeta_{\alpha,1})^{n+1}}t^{n}\\
&\,-\frac{E_{q}(-2i\zeta_{\alpha,1}z)E_{q}(i\zeta_{\alpha,1})}
{(1-q)\frac{d}{dz}\mathcal{J}^{(2)}_{\alpha}(z;q^{2})|_{z=j^{(2)}_{1,\alpha}}}\sum_{n=0}^{\infty}\frac{(-i)^{-n}}{(2\zeta_{\alpha,1})^{n+1}}t^{n}.
\end{split}
\end{equation}

Substituting with $ \pm i=e^{\pm\frac{i\pi}{2}}$ into (\ref{y787}), we obtain
\begin{equation}\label{s8}
\begin{split}
  g(z,t)=-\frac{2}{1-q}\,Re\left(\frac{E_{q}(2i\zeta_{\alpha,1}z)E_{q}(-i\zeta_{\alpha,1})}{\frac{d}{dz}\mathcal{J}^{(2)}_{\alpha}(z;q^{2})|_{z=j^{(2)}_{1,\alpha}}}
\sum_{n=0}^{\infty}\frac{e^{\frac{-in\pi}{2}}}{(2\zeta_{\alpha,1})^{n+1}}t^{n}\right).
\end{split}
\end{equation}
But

\begin{equation}\label{s14}
\begin{split}
Re\,\left( e^{\frac{-in\pi}{2}}E_{q}(2i\zeta_{\alpha,1}z) E_{q}(-i\zeta_{\alpha,1})\right)&=\cos\frac{n\pi}{2}
\left(Cos_{q}(2\zeta_{\alpha,1}z)Cos_{q}(\zeta_{\alpha,1})+Sin_{q}(2\zeta_{\alpha,1}z)Sin_{q}(\zeta_{\alpha,1})\right)\\&+
\sin\frac{n\pi}{2}\left(Sin_{q}(2\zeta_{\alpha,1}z)Cos_{q}(\zeta_{\alpha,1})
-Cos_{q}(2\zeta_{\alpha,1}z)Sin_{q}(\zeta_{\alpha,1})\right).
 \end{split}
\end{equation}

Therefore, substituting  from (\ref{s14}) into (\ref{s8}), we get for $|t|=2\zeta_{\alpha,1}$

\begin{equation*}
\begin{split}
  g(z,t)=&-\frac{2}{1-q}\frac{
Cos_{q}(2\zeta_{\alpha,1}z)Cos_{q}(\zeta_{\alpha,1})+Sin_{q}(2\zeta_{\alpha,1}z)Sin_{q}(\zeta_{\alpha,1})}
{\frac{d}{dz}\mathcal{J}^{(2)}_{\alpha}(z;q^{2})|_{z=j^{(2)}_{1,\alpha}}}
\sum_{n=0}^{\infty}\frac{\,\cos\frac{n\pi}{2}}{(2\zeta_{\alpha,1})^{n+1}}t^{n}\\&-\frac{2}{1-q}
\frac{Sin_{q}(2\zeta_{\alpha,1}z)Cos_{q}(\zeta_{\alpha,1})
-Cos_{q}(2\zeta_{\alpha,1}z)Sin_{q}(\zeta_{\alpha,1})}{\frac{d}{dz}\mathcal{J}^{(2)}_{\alpha}(z;q^{2})|_{z=j^{(2)}_{1,\alpha}}}
\sum_{n=0}^{\infty}\frac{\,\sin\frac{n\pi}{2}}{(2\zeta_{\alpha,1})^{n+1}}t^{n}.
\end{split}
\end{equation*}

Consequently,

\begin{equation*}
\begin{split}
  b_{n}=&-\frac{2}{1-q}\cos\frac{n\pi}{2}\frac{
\left(Cos_{q}(2\zeta_{\alpha,1}z)Cos_{q}(\zeta_{\alpha,1})+Sin_{q}(2\zeta_{\alpha,1}z)Sin_{q}(\zeta_{\alpha,1})\right)}{(2\zeta_{\alpha,1})^{n+1}
\frac{d}{dz}\mathcal{J}^{(2)}_{\alpha}(z;q^{2})|_{z=j^{(2)}_{1,\alpha}}
}\\&-\frac{2}{1-q}\sin\frac{n\pi}{2}
\frac{\left(Sin_{q}(2\zeta_{\alpha,1}z)Cos_{q}(\zeta_{\alpha,1})
-Cos_{q}(2\zeta_{\alpha,1}z)Sin_{q}(\zeta_{\alpha,1})\right)}{(2\zeta_{\alpha,1})^{n+1}\frac{d}{dz}\mathcal{J}^{(2)}_{\alpha}(z;q^{2})|_{z=j^{(2)}_{1,\alpha}}
}.\end{split}
\end{equation*}

By applying the Darboux method, we obtain
 $$a_{n}=b_{n}+o(r^{-n})\quad(n\rightarrow\infty).$$
Thus

\begin{equation*}
\begin{split}
\frac{ B^{(2)}_{n,\alpha}(z;q)}{[n]_{q}!}
=&-\frac{2}{1-q}\cos\frac{n\pi}{2}\frac{\left(Cos_{q}(2\zeta_{\alpha,1}z)Cos_{q}(\zeta_{\alpha,1})+Sin_{q}(2\zeta_{\alpha,1}z)Sin_{q}(\zeta_{\alpha,1})\right)}
{(2\zeta_{\alpha,1})^{n+1}\frac{d}{dz}\mathcal{J}^{(2)}_{\alpha}(z;q^{2})|_{z=j^{(2)}_{1,\alpha}}}
\\&-\frac{2}{1-q}\sin\frac{n\pi}{2}\frac{\left(Sin_{q}(2\zeta_{\alpha,1}z)Cos_{q}(\zeta_{\alpha,1})
-Cos_{q}(2\zeta_{\alpha,1}z)Sin_{q}(\zeta_{\alpha,1})\right)}{(2\zeta_{\alpha,1})^{n+1}
\frac{d}{dz}\mathcal{J}^{(2)}_{\alpha}(z;q^{2})|_{z=j^{(2)}_{1,\alpha}}}+o(r^{-n}).\end{split}
\end{equation*}

Hence

\begin{align*}
\begin{split}
B^{(2)}_{2n,\alpha}(z;q)=&\frac{2}{1-q}\frac{(-1)^{n+1}{[2n]_{q}!}\left(Cos_{q}(2\zeta_{\alpha,1}z)Cos_{q}(\zeta_{\alpha,1})+Sin_{q}(2\zeta_{\alpha,1}z)Sin_{q}(\zeta_{\alpha,1})
\right)}{(2\zeta_{\alpha,1})^{2n+1}\frac{d}{dz}\mathcal{J}^{(2)}_{\alpha}(z;q^{2})|_{z=j^{(2)}_{1,\alpha}}}\\&+o(r^{-2n})
,
 \\
 B^{(2)}_{2n+1,\alpha}(z;q)=&\frac{2}{1-q}\frac{(-1)^{n+1}{[2n+1]_{q}!}\left(Sin_{q}(2\zeta_{\alpha,1}z)Cos_{q}(\zeta_{\alpha,1})
-Cos_{q}(2\zeta_{\alpha,1}z)Sin_{q}(\zeta_{\alpha,1})\right)}{(2\zeta_{\alpha,1})^{2n+2}\frac{d}{dz}\mathcal{J}^{(2)}_{\alpha}(z;q^{2})|_{z=j^{(2)}_{1,\alpha}}}
\\&+o(r^{-2n}),\\
 \end{split}
\end{align*}
which  are the required results.
\end{proof}

\begin{corollary}
Let  $(\beta_{n,\alpha}(q))_{n}$  be the generalized $q$-Bernoulli numbers   which are  defined in (\ref{q66}). Then
\begin{align*}
\begin{split}
\beta_{2n,\alpha}(q)&=\dfrac{2}{1-q}\dfrac{(-1)^{n+1}[2n]_{q}!
Cos_{q}(\zeta_{\alpha,1})
}{(2\zeta_{\alpha,1})^{2n+1}\,\frac{d}{dz}\mathcal{J}^{(2)}_{\alpha}(z;q^{2})|_{z=j^{(2)}_{1,\alpha}}
}+o(r^{-2n}),\\
\beta_{2n+1,\alpha}(q)&=\dfrac{2}{1-q}\dfrac{(-1)^{n}[2n+1]_{q}!
Sin_{q}(\zeta_{\alpha,1})}{(2\zeta_{\alpha,1})^{2n+2}\,\frac{d}{dz}\mathcal{J}^{(2)}_{\alpha}(z;q^{2})|_{z=j^{(2)}_{1,\alpha}}}+o(r^{-2n})
,\end{split}
\end{align*}
as  $n\,\rightarrow\,\infty,$ where $0<r\leq 2\zeta_{\alpha,1}$  and  $\zeta_{\alpha,1}=\dfrac{j^{(2)}_{1,\alpha}}{2(1-q)}$.
\end{corollary}
\begin{proof}
The proof follows from Theorem \ref{rr} by substituting with $z=0$.
\end{proof}
\vskip0.1cm
\begin{corollary}
Let  $(B_{n}(z;q))_{n}$  be  the $q$-Bernoulli polynomials  which are  defined in (\ref{q8}).  Then
\begin{align*}
B_{2n}(z;q)&=\dfrac{1}{1-q}\dfrac{(-1)^{n+1}[2n]_{q}!
Cos_{q}(2\zeta_{1,q}z)Cos_{q}(\zeta_{1,q})}{(2\zeta_{1,q})^{2n}Sin^{'}_{q}(\zeta_{1,q})
}+o(r^{-2n}),\\
B_{2n+1}(z;q)&=\dfrac{1}{1-q}\dfrac{(-1)^{n+1}[2n+1]_{q}!Sin_{q}(2\zeta_{1,q}z)Cos_{q}(\zeta_{1,q})
}{(2\zeta_{1,q})^{2n+1}Sin^{'}_{q}(\zeta_{1,q})}+o(r^{-2n})
,
\end{align*}
as  $n\,\rightarrow\,\infty$, where $0<r\leq 2\zeta_{1,q}$  and  $\zeta_{1,q}$ is the smallest positive zero of $Sin_{q}z$.

\end{corollary}
\vskip0.1cm
\begin{proof}
The proof follows by setting $\alpha=\dfrac{1}{2}$ in Theorem \ref{rr} and using that $\zeta_{1,q}=\dfrac{j^{(2)}_{1,1/2}}{2(1-q)}$ is the smallest postive zero of $Sin_{q}z $.
\end{proof}
\vskip0.1cm
\begin{corollary}
Let $(E_{n}(z;q))_{n }$  be  the $q$-Euler polynomials  which are  defined in (\ref{q13248}). Then
\begin{align*}
E_{2n}(z;q)&=\dfrac{2}{1-q}\dfrac{(-1)^{n+1}[2n]_{q}!
Sin_{q}(2\eta_{1,q}z)Sin_{q}(\eta_{1,q})
}{(2\eta_{1,q})^{2n+1}Cos^{'}_{q}(\eta_{1,q})
}+o(r^{-2n}),\\
E_{2n+1}(z;q)&=\dfrac{2}{1-q}\dfrac{(-1)^{n}[2n+1]_{q}!
Cos_{q}(2\eta_{1,q}z)Sin_{q}(\eta_{1,q})}{(2\eta_{1,q})^{2n+2}Cos^{'}_{q}(\eta_{1,q})}+o(r^{-2n})
,
\end{align*}
as  $n\,\rightarrow\,\infty,$  where $0<r\leq 2\eta_{1,q}$ and  $\eta_{1,q}$ is  the smallest positive  zero of $Cos_{q}z$.

\end{corollary}
\vskip0.1cm
\begin{proof}
The proof follows by substituting with $\alpha=\dfrac{-1}{2}$ in Theorem \ref{rr} and using that $\eta_{1,q}=\dfrac{j^{(2)}_{1,-1/2}}{2(1-q)}$ is a zero of $Cos_{q}z $.
\end{proof}
\vskip0.1cm

\begin{theorem}\label{mm}
 The large $n$ asymptotic of the generalized $q$-Bernoulli polynomials ${\small(B^{(3)}_{n,\alpha}(z;q))_{n}}$ are
\begin{equation*}
\begin{split}
B^{(3)}_{2n,\alpha}(z;q)&=\dfrac{4q^{\frac{1}{4}}}{1-q}\dfrac{(-1)^{n+1}[2n]_{q}!
\left(C_{q}(2\eta_{\alpha,1}z)C_{q}(\eta_{\alpha,1})+S_{q}(2\eta_{\alpha,1}z)S_{q}(\eta_{\alpha,1})\right)
}{(2\eta_{\alpha,1})^{2n+1}\,\frac{d}{dz}\mathcal{J}^{(3)}_{\alpha}(z;q^{2})|_{z=j^{(3)}_{1,\alpha}}
}\\&+o(r^{-2n}),
\end{split}
\end{equation*}
\begin{equation*}
\begin{split}
\quad B^{(3)}_{2n+1,\alpha}(z;q)&=\dfrac{4q^{\frac{1}{4}}}{1-q}\dfrac{(-1)^{n+1}[2n+1]_{q}!
\left(S_{q}(2\eta_{\alpha,1}z)C_{q}(\eta_{\alpha,1})
-C_{q}(2\eta_{\alpha,1}z)S_{q}(\eta_{\alpha,1})\right)
}{(2\eta_{\alpha,1})^{2n+2}\,\frac{d}{dz}\mathcal{J}^{(3)}_{\alpha}(z;q^{2})|_{z=j^{(3)}_{1,\alpha}}}\\&+o(r^{-2n})
,
\end{split}
\end{equation*}
as  $n\,\rightarrow\,\infty,$ where $S_{q}(\cdot)$ and $C_{q}(\cdot)$  are the functions defined in (\ref{h765}), $0<r\leq 2\eta_{\alpha,1},$ and  $\eta_{\alpha,1}:=\dfrac{q^{\frac{1}{4}}j^{(3)}_{1,\alpha}}{1-q}$.
\end{theorem}
\vskip0.1cm
\begin{proof}
 Since the generalized $q$-Bernoulli polynomials $B^{(3)}_{n,\alpha}(z;q)$ are defined by the generating function
\begin{equation*}
  f(z,t):=\frac{\exp_{q}(zt)\exp_{q}(\frac{-t}{2})}{g^{(3)}_{\alpha}(it;q)}=\sum_{n=0}^{\infty}B^{(3)}_{n,\alpha}(z;q)\frac{t^{n}}{[n]_{q}!},
\end{equation*}
 where $f(z,t)$ is analytic in the disk  $|t|<\frac{2q^{\frac{1}{4}}j^{(3)}_{1,\alpha}}{1-q}$. We take
\begin{equation*}
\begin{split}
  g(z,t)&:=Res(f,\frac{2iq^{\frac{1}{4}}j^{(3)}_{1,\alpha}}{1-q})\frac{1}{(t-\frac {2iq^{\frac{1}{4}}j^{(3)}_{1,\alpha}}{1-q})}+Res(f,\frac{-2iq^{\frac{1}{4}}j^{(3)}_{1,\alpha}}{1-q})
  \frac{1}{(t+\frac {2iq^{\frac{1}{4}}j^{(3)}_{1,\alpha}}{1-q})}.
\end{split}
\end{equation*}
Set
\begin{equation}\label{h763}
  f(z,t)= \displaystyle\sum_{n=0}^{\infty}a_{n}\, t^{n}, \quad \text {and} \quad   g(z,t)=\displaystyle\sum_{n=0}^{\infty}b_{n}\, t^{n}.
\end{equation}
 Then
\begin{equation*}
\begin{split}
  g(z,t)=Res(f,2i\eta_{\alpha,1})\frac{1}{(t-2i\eta_{\alpha,1})} +Res(f,-2i\eta_{\alpha,1})\frac{1}{(t+2i\eta_{\alpha,1})}.
\end{split}
\end{equation*}

 The function  $f(z,t)-g(z,t)$ is continuous within and on the the circle $|t|=2\eta_{\alpha,1}$.  Now, we calculate the coefficient of $t^{n}$ in the expansion of $g(z,t)$. Indeed,
\begin{equation}\label{y70}
\begin{split}
  g(z,t)&=Res(f,2i\eta_{\alpha,1})\frac{1}{(t-2i\eta_{\alpha,1})} +Res(f,-2i\eta_{\alpha,1})\frac{1}{(t+2i\eta_{\alpha,1})}
\\& =Res(f,2i\eta_{\alpha,1})\frac{-1}{2i\eta_{\alpha,1}}\sum_{n=0}^{\infty}(\frac{t}{2i\eta_{\alpha,1}})^{n}+Res(f,-2i\eta_{\alpha,1})\frac{1}{2i\eta_{\alpha,1}}
\sum_{n=0}^{\infty}(\frac{-t}{2i\eta_{\alpha,1}})^{n}.
\end{split}
\end{equation}
Since

\begin{equation*}
\begin{split}
&Res(f,2i\eta_{\alpha,1} )=\frac{\exp_{q}(2i\eta_{\alpha,1}z)\exp_{q}(-i\eta_{\alpha,1})}{\frac{d}{dz}\mathcal{J}^{(3)}_{\alpha}(\frac{iz}{2}(1-q)q^{\frac{-1}{4}};q^{2})|_
{z=2i\eta_{\alpha,1}}},\\&
Res(f,-2i\eta_{\alpha,1} )=\frac{\exp_{q}(-2i\eta_{\alpha,1}z)\exp_{q}(i\eta_{\alpha,1})}
{\frac{d}{dz}\mathcal{J}^{(3)}_{\alpha}(\frac{iz}{2}(1-q)q^{\frac{-1}{4}};q^{2})|_{z=-2i\eta_{\alpha,1}}},
 \end{split}
\end{equation*}
then Equation (\ref{y70}) can be written as
\begin{equation*}
\begin{split}
  g(z,t)=&-\frac{\exp_{q}(2i\eta_{\alpha,1}z)\exp_{q}(-i\eta_{\alpha,1})}{\frac{d}{dz}\mathcal{J}^{(3)}_{\alpha}(\frac{iz}{2}(1-q)q^{\frac{-1}{4}};q^{2})|_{z=2i\eta_{\alpha,1}}}
\sum_{n=0}^{\infty}\frac{i^{-n-1}}{(2\eta_{\alpha,1})^{n+1}}t^{n}
\\&-\frac{\exp_{q}(-2i\eta_{\alpha,1}z)\exp_{q}(i\eta_{\alpha,1})}
{\frac{d}{dz}\mathcal{J}^{(3)}_{\alpha}(\frac{iz}{2}(1-q)q^{\frac{-1}{4}};q^{2})|_{z=-2i\eta_{\alpha,1}}}\sum_{n=0}^{\infty}\frac{(-i)^{-n-1}}
{(2\eta_{\alpha,1})^{n+1}}t^{n}.
\end{split}
\end{equation*}
Since
\begin{equation*}
\begin{split}
\frac{d}{dz}\mathcal{J}^{(3)}_{\alpha}(\frac{iz}{2}(1-q)q^{\frac{-1}{4}};q^{2})|_{z=2i\eta_{\alpha,1}}&=\frac{i}{2}(1-q)q^{\frac{-1}{4}}\frac{d}{dz}\mathcal{J}^{(3)}_{\alpha}(\frac{iz}{2}(1-q)
q^{\frac{-1}{4}};q^{2})|_{z=2i\eta_{\alpha,1}}\\&=-\frac{i}{2}(1-q)
q^{\frac{-1}{4}}\frac{d}{dz}\mathcal{J}^{(3)}_{\alpha}(z;q^{2})|_{z=j^{(3)}_{\alpha,1}},\end{split}
\end{equation*}
then
\begin{equation}\label{y7087}
\begin{split}
  g(z,t)=&-\frac{2q^{\frac{-1}{4}}}{1-q}\frac{\exp_{q}(2i\eta_{\alpha,1}z)\exp_{q}(-i\eta_{\alpha,1})}{\frac{d}{dz}\mathcal{J}^{(3)}_{\alpha}(z;q^{2})
  |_{z=j^{(3)}_{\alpha,1}}}
\sum_{n=0}^{\infty}\frac{i^{-n}}{(2\eta_{\alpha,1})^{n+1}}t^{n}
\\&-\frac{2q^{\frac{-1}{4}}}{1-q}\frac{\exp_{q}(-2i\eta_{\alpha,1}z)\exp_{q}(i\eta_{\alpha,1})}
{\frac{d}{dz}\mathcal{J}^{(3)}_{\alpha}(z;q^{2})|_{z=j^{(3)}_{\alpha,1}}}\sum_{n=0}^{\infty}\frac{(-i)^{-n}}{(2\eta_{\alpha,1})^{n+1}}t^{n}.
\end{split}
\end{equation}
Substituting with $ \pm i=e^{\pm\frac{i\pi}{2}}$ into (\ref{y7087}),  we get
 \begin{equation*}
\begin{split}
  g(z,t)=&-\frac{2q^{\frac{1}{4}}}{1-q}\frac{\exp_{q}(2i\eta_{\alpha,1}z)\exp_{q}(-i\eta_{\alpha,1})}
  {\frac{d}{dz}\mathcal{J}^{(3)}_{\alpha}(z;q^{2})|_{z=j^{(3)}_{1,\alpha}}}
\sum_{n=0}^{\infty}\frac{e^{\frac{-in\pi}{2}}}{(2\eta_{\alpha,1})^{n+1}}t^{n}\\&-\frac{2q^{\frac{1}{4}}}{1-q}
\frac{\exp_{q}(-2i\eta_{\alpha,1}z)\exp_{q}(i\eta_{\alpha,1})}
{\frac{d}{dz}\mathcal{J}^{(3)}_{\alpha}(z;q^{2})|_{z=j^{(3)}_{1,\alpha}}}\sum_{n=0}^{\infty}\frac{e^{\frac{in\pi}{2}}}
{(2\eta_{\alpha,1})^{n+1}}t^{n}.
\end{split}
\end{equation*}
This implies
\begin{equation}\label{s98}
\begin{split}
  g(z,t)=\frac{-4q^{\frac{1}{4}}}{1-q}\,Re\left(\frac{\exp_{q}(2i\eta_{\alpha,1}z)\exp_{q}(-i\eta_{\alpha,1})}{\frac{d}{dz}\mathcal{J}^{(3)}_{\alpha}(z;q^{2})|_{z=j^{(3)}_{1,\alpha}}}
\sum_{n=0}^{\infty}\frac{e^{\frac{-in\pi}{2}}}{(2\eta_{\alpha,1})^{n+1}}t^{n}\right).
\end{split}
\end{equation}
But
\begin{equation}\label{s124}
\begin{split}
Re\,\left(e^{\frac{-in\pi}{2}} \exp_{q}(2i\eta_{\alpha,1}z) \exp_{q}(-i\eta_{\alpha,1})\right)&=\cos\frac{n\pi}{2}
\left(C_{q}(2\eta_{\alpha,1}z)C_{q}(\eta_{\alpha,1})+S_{q}(2\eta_{\alpha,1}z)S_{q}(\eta_{\alpha,1})\right)\\&+
\sin\frac{n\pi}{2}\left(S_{q}(2\eta_{\alpha,1}z)C_{q}(\eta_{\alpha,1})
-C_{q}(2\eta_{\alpha,1}z)S_{q}(\eta_{\alpha,1})\right).
 \end{split}
\end{equation}
Therefore, substituting  from (\ref{s124}) into (\ref{s98}), we get for $|t|=2\eta_{\alpha,1}$
\begin{equation*}
\begin{split}
  g(z,t)=&-\frac{4q^{\frac{1}{4}}}{1-q}\frac{
\left(C_{q}(2\eta_{\alpha,1}z)C_{q}(\eta_{\alpha,1})+S_{q}(2\eta_{\alpha,1}z)S_{q}(\eta_{\alpha,1})\right)}{\,
\frac{d}{dz}\mathcal{J}^{(3)}_{\alpha}(z;q^{2})|_{z=j^{(3)}_{1,\alpha}}
}\sum_{n=0}^{\infty}\frac{\,\cos\frac{n\pi}{2}}{(2\eta_{\alpha,1})^{n+1}}t^{n}\\&-\frac{4q^{\frac{1}{4}}}{1-q}
\frac{\left(S_{q}(2\eta_{\alpha,1}z)C_{q}(\eta_{\alpha,1})
-C_{q}(2\eta_{\alpha,1}z)S_{q}(\eta_{\alpha,1})\right)}{\frac{d}{dz}\mathcal{J}^{(3)}_{\alpha}(z;q^{2})|_{z=j^{(3)}_{1,\alpha}}
}\sum_{n=0}^{\infty}\frac{\,\sin\frac{n\pi}{2}}{(2\eta_{\alpha,1})^{n+1}}t^{n}.
\end{split}
\end{equation*}
Hence
\begin{equation*}
\begin{split}
  b_{n}=&-\frac{4q^{\frac{1}{4}}}{1-q}\cos\frac{n\pi}{2}\frac{
\left(C_{q}(2\eta_{\alpha,1}z)C_{q}(\eta_{\alpha,1})+S_{q}(2\eta_{\alpha,1}z)S_{q}(\eta_{\alpha,1})\right)}{(2\eta_{\alpha,1})^{n+1}
\,\frac{d}{dz}\mathcal{J}^{(3)}_{\alpha}(z;q^{2})|_{z=j^{(3)}_{1,\alpha}}}\\&-\frac{4q^{\frac{1}{4}}}{1-q}\sin\frac{n\pi}{2}
\frac{\left(S_{q}(2\eta_{\alpha,1}z)C_{q}(\eta_{\alpha,1})
-C_{q}(2\eta_{\alpha,1}z)S_{q}(\eta_{\alpha,1})\right)}{(2\eta_{\alpha,1})^{n+1}\,\frac{d}{dz}\mathcal{J}^{(3)}_{\alpha}(z;q^{2})|_
{z=j^{(3)}_{1,\alpha}}
}.\end{split}
\end{equation*}
By applying the Darboux method, we get
 $$a_{n}=b_{n}+o(r^{-n})\quad(n\rightarrow\infty).$$
Thus
\begin{equation*}
\begin{split}
\frac{ B^{(3)}_{n,\alpha}(z;q)}{[n]_{q}!}
=&-\frac{4q^{\frac{1}{4}}}{1-q}\cos\frac{n\pi}{2}\frac{
\left(C_{q}(2\eta_{\alpha,1}z)C_{q}(\eta_{\alpha,1})+S_{q}(2\eta_{\alpha,1}z)S_{q}(\eta_{\alpha,1})\right)}{(2\eta_{\alpha,1})^{n+1}
\,\frac{d}{dz}\mathcal{J}^{(3)}_{\alpha}(z;q^{2})|_{z=j^{(3)}_{1,\alpha}}}\\&-\frac{4q^{\frac{1}{4}}}{1-q}\sin\frac{n\pi}{2}
\frac{\left(S_{q}(2\eta_{\alpha,1}z)C_{q}(\eta_{\alpha,1})
-C_{q}(2\eta_{\alpha,1}z)S_{q}(\eta_{\alpha,1})\right)}{(2\eta_{\alpha,1})^{n+1}\,
\frac{d}{dz}\mathcal{J}^{(3)}_{\alpha}(z;q^{2})|_{z=j^{(3)}_{1,\alpha}}
}+o(r^{-n}).\end{split}
\end{equation*}
Therefore,
\begin{align*}\label{s1432}
\begin{split}
B^{(3)}_{2n,\alpha}(z;q)&=\frac{4q^{\frac{1}{4}}}{1-q}\frac{(-1)^{n+1}[2n]_{q}!\left(C_{q}(2\eta_{\alpha,1}z)C_{q}(\eta_{\alpha,1})+S_{q}(2\eta_{\alpha,1}z)S_{q}(\eta_{\alpha,1})\right)}
{(2\eta_{\alpha,1})^{2n+1}\,\frac{d}{dz}\mathcal{J}^{(3)}_{\alpha}(z;q^{2})|_{z=j^{(3)}_{1,\alpha}}}\\&+o(r^{-2n})
,
 \\
 B^{(3)}_{2n+1,\alpha}(z;q)&=\frac{4q^{\frac{1}{4}}}{1-q}\frac{(-1)^{n+1}[2n+1]_{q}!\left(S_{q}(2\eta_{\alpha,1}z)C_{q}(\eta_{\alpha,1})
-C_{q}(2\eta_{\alpha,1}z)S_{q}(\eta_{\alpha,1})\right)}{(2\eta_{\alpha,1})^{2n+2}\,
\frac{d}{dz}\mathcal{J}^{(3)}_{\alpha}(z;q^{2})|_{z=j^{(3)}_{1,\alpha}}}\\&+o(r^{-2n})
,
 \end{split}
\end{align*}
which are the desired result.

\end{proof}
\vskip0.1cm
\begin{corollary}
The large $n$ asymptotic of the generalized $q$-Bernoulli numbers $(\beta^{(3)}_{n,\alpha}(q))_{n}$ defined in (\ref{r66}) are
\begin{equation*}
\begin{split}
\beta^{(3)}_{2n,\alpha}(q)=\dfrac{4q^{\frac{1}{4}}}{1-q}\dfrac{(-1)^{n+1}[2n]_{q}!
C_{q}(\eta_{\alpha,1})}{(2\eta_{\alpha,1})^{2n+1}\,\frac{d}{dz}\mathcal{J}^{(3)}_{\alpha}(z;q^{2})|_{z=j^{(3)}_{1,\alpha}}
}+o(r^{-2n}),
\end{split}
\end{equation*}
\begin{equation*}
\begin{split}
\beta^{(3)}_{2n+1,\alpha}(q)=\dfrac{4q^{\frac{1}{4}}}{1-q}\dfrac{(-1)^{n}[2n+1]_{q}!
S_{q}(\eta_{\alpha,1})}{(2\eta_{\alpha,1})^{2n+2}\,\frac{d}{dz}\mathcal{J}^{(3)}_{\alpha}(z;q^{2})|_{z=j^{(3)}_{1,\alpha}}}+o(r^{-2n})
,
\end{split}
\end{equation*}
as  $n\,\rightarrow\,\infty$,  where  $0<r\leq 2\eta_{\alpha,1}$  and   $\eta_{\alpha,1}=\dfrac{q^{\frac{1}{4}}j^{(3)}_{1,\alpha}}{1-q}$.

\end{corollary}
\vskip0.1cm
\begin{proof}
The proof follows from Theorem \ref{mm} by substituting with $z=0$.
\end{proof}
\vskip0.1cm
\begin{corollary}
The large $n$ asymptotic of $q$-Bernoulli polynomials $(\tilde{B}_{n}(z;q))_{n}$  defined in (\ref{q4039}) are
\begin{align*}
\tilde{B}_{2n}(z;q)&=\dfrac{2q^{\frac{1}{4}}}{1-q}\dfrac{(-1)^{n+1}[2n]_{q}!
C_{q}(2\lambda_{1,q}z)C_{q}(\lambda_{1,q})}{(2\lambda_{1,q})^{2n}(S_{q})^{'}(\lambda_{1,q})
}+o(r^{-2n}),\\
\tilde{B}_{2n+1}(z;q)&=\dfrac{2q^{\frac{1}{4}}}{1-q}\dfrac{(-1)^{n+1}[2n+1]_{q}!S_{q}(2\lambda_{1,q}z)C_{q}(\lambda_{1,q})
}{(2\lambda_{1,q})^{2n+1}(S_{q})^{'}(\lambda_{1,q})}+o(r^{-2n})
,
\end{align*}
as  $n\,\rightarrow\,\infty$, where $0<r\leq 2\lambda_{1,q}$ and  $\lambda_{1,q}$ is the smallest positive zero of $S_{q}(z)$.

\end{corollary}
\vskip0.1cm
\begin{proof}
The proof follows by substituting with  $\alpha=\dfrac{1}{2}$ in Theorem \ref{mm} and using that $\lambda_{1,q}=\dfrac{q^{\frac{1}{4}}j^{(3)}_{1,1/2}}{1-q}$ is a zero of $S_{q}(z)$.
\end{proof}
\vskip0.1cm
\begin{corollary}
The large $n$ asymptotic of the $q$-Euler polynomials $(\tilde{E}_{n}(z;q))_{n }$  defined in (\ref{q4039}) are
\begin{align*}
\tilde{E}_{2n}(z;q)&=\dfrac{4q^{\frac{1}{4}}}{1-q}\dfrac{(-1)^{n+1}[2n]_{q}!
S_{q}(2\mu_{1,q}z)S_{q}(\mu_{1,q})
}{(2\mu_{1,q})^{2n+1}(C_{q})^{'}(q^{\frac{1}{2}}\mu_{1,q})
}+o(r^{-2n}),\\
\tilde{E}_{2n+1}(z;q)&=\dfrac{4q^{\frac{1}{4}}}{1-q}\dfrac{(-1)^{n}[2n+1]_{q}!
C_{q}(2\mu_{1,q}z)S_{q}(\mu_{1,q})}{(2\mu_{1,q})^{2n+2}(C_{q})^{'}(q^{\frac {1}{2}}\mu_{1,q})}+o(r^{-2n})
,
\end{align*}
as  $n\,\rightarrow\,\infty,$  where  $0<r\leq 2\mu_{1,q}$   and  $\mu_{1,q}$ is the smallest positive zero of $C_{q}(q^{\frac{1}{2}}z)$.

\end{corollary}
\vskip0.1cm
\begin{proof}
The proof follows by setting $\alpha=\dfrac{-1}{2}$ in Theorem \ref{mm} and using that $\mu_{1,q}=\dfrac{q^{\frac{1}{4}}j^{(3)}_{1,-1/2}}{1-q}$ is a zero of $C_{q}(q^{\frac{1}{2}}z)$.
\end{proof}
\section{Expansion of the generalized $q$-Bernoulli polynomials}
 This section introduces an expansion  for a class of  entire function $\emph{f}$  in terms of the generalized $q$-Bernoulli polynomials defined in (\ref{q14})-(\ref{q197}). We need the following definitions and results.
\begin{definition}
 Let $f$ be an entire function and $M(r;f):=\max\{|f(z)|:|z|=r\}$.  Then the order $\rho$ of $f$ is defined by
\begin{equation*}
  \rho=\limsup_{r\rightarrow\infty}\frac{\log(\log M(r;f))}{\log r}.
\end{equation*}
If $0<\rho < \infty$, the type $\sigma$ is defined as
\begin{equation*}
 \sigma=\limsup_{r\rightarrow\infty}\frac{\log M(r;f)}{r^{\rho}},
\end{equation*}
see \cite{boass1}.
\end{definition}
\begin{definition}
 Let $k$ be a non zero real number, and let $p$ be a real
number with $|p| > 1$. An entire function $f$ has a $p$-exponential growth
of order $k$ and a finite type $\gamma$, if there exists real numbers $K > 0$ and $\gamma$, such that
\begin{equation*}
  |f(z)| \leq Kp^{\frac{k}{2}(\frac{log|z|}{log p})^{2}}|z|^{\gamma},
\end{equation*}
or equivalently,
\begin{equation*}
  |f(z)|\leq Ke^{\frac{k}{2 log p}(log |z|)^{2}+\gamma log|z|},
\end{equation*}
see \cite{ramis}.
\end{definition}
\vskip0.1cm
 In \cite[Lemma 2.2]{ramis}, Ramis proved that an entire function $ f(z)=\displaystyle\sum_{n=0}^{\infty}a_{n}z^{n}$ has a $p$-exponential growth of order $k$ and a finite type $\gamma$ if and only if there exists a constant $K>0$ such that
 \begin{equation}\label{r65}
  |a_{n}|<K \,p^{\frac{-(n-\gamma)^{2}}{2k}}\quad (n\in\mathbb{N}_{0}).
\end{equation}
\begin{definition}
The function $f$ is of $\psi$-type $\tau$  if there exists constants $M$ and $\tau$ such that
\begin{equation}
 f(r e^{i\theta})\leq M \Psi(\tau r)\quad (r\rightarrow\infty).
\end{equation}
\end{definition}
\vskip0.1cm

In the following, we assume that $\Psi$ is a comparison function, i.e.
\[\Psi(t)=\displaystyle\sum_{n=0}^{\infty} \psi_{n} t^{n},\,\,\psi_{n}\neq 0,\;\mbox{for all}\;
 n\in\mathbb{N}_0,\]
 and  \[\displaystyle\lim_{n\rightarrow\infty}\frac{\Psi_{n+1}}{\Psi_{n}}=0,\] see \cite{boass}.  Nachbin in \cite{nac}  proved that a function $ f(z)=\displaystyle\sum_{n=0}^{\infty}a_{n}z^{n}$   of $\psi$-type $\tau$ if and only if
\begin{equation}\label{r54}
  \displaystyle\limsup_{n\rightarrow\infty}\sqrt[n]{|\frac{f_{n}}{\psi_{n}}|}=\tau.
\end{equation}

The Borel transform of a function $f(z)=\displaystyle\sum_{n=0}^{\infty}f_{n}z^{n}$  is defined by
\begin{equation}\label{r32}
 F(w)=\displaystyle\sum_{n=0}^{\infty} \frac{f_{n}}{\Psi_{n}}\frac{1}{w^{n+1}} .
\end{equation}
 Hence $F(w)$ is analytic if $\displaystyle\lim_{n\rightarrow\infty}\sqrt[n]{|\frac{f_{n}}{\Psi_{n}}|}\frac{1}{|w|}<1$,  i.e $F(w)$ is analytic for $|w|>\tau$.
\vskip0.1cm
 In \cite{boass}, a set of polynomials $\{P_{n}(z)\}_{n}$ has a generalized Appell representation if it is generated by the formal relation

 \begin{equation}\label{w2}
  B(\xi)\Psi(z\xi)=\sum_{n=0}^{\infty}P_{n}(z)( W(\xi))^{n}.
\end{equation}

Here $B(\xi)$ and $W(\xi)$ are analytic at zero. We choose $\Omega$ to be a region in which $B(\xi)$ and $W(\xi)$ are analytic. If $\rho_{0}$ is the distance from the origin to the nearest of the boundary of $\Omega$, the series in (\ref{w2}) is uniformly convergent in compact subset of $\Delta$, where $\Delta$ is the open disk $\{\xi\in\Omega:|\xi|<\rho_{0}\}$.

\vskip 0.2 cm
In the following, let $C$  be a compact subset of $\Delta$ which include $D(f)$, where $D(f)$ is the union of the set of all singular points of $F$ and the set of all points exterior to the domain of $F$. Let $\Re_{\Psi}[C]$  be the class of entire functions of finite $\Psi$-type less than or equal $\tau$.

\vskip 0.2 cm

The following theorem  from \cite[Theorem 7.2, P. 22]{boass} is essential in our investigation.

\begin{theorem}\label{jj}
Let $(P_{n}(z))_{n}$ be the generalized Appell polynomials defined in (\ref{w2}). If $f\in\Re_{\Psi}[C]$,  then
\begin{equation*}
  f(z)=\sum_{n=0}^{\infty} L_{n}(f)P_{n}(z),
\end{equation*}
where
\begin{equation*}
   L_{n}(f)=\frac{1}{2\pi i}\int_{\Gamma}\frac{w^{n}}{B(w)}F(w)dw,
\end{equation*}
 and  $F(w)$ is the Borel transform of $f(z)$  defined in (\ref{r32}) and $\Gamma$ is a circle $|w|=\rho$ with $\tau<\rho<\rho_{0}$.
\end{theorem}

\vskip 0.1 cm

We derive the following  theorem by applying Theorem \ref{jj} on the polynomials $(B^{(2)}_{n,\alpha}(\cdot;q))_{n}$ defined in (\ref{q13}).
\begin{theorem}\label{bb}
Let  $f(z)=\displaystyle\sum_{n=0}^{\infty}f_{n}z^{n}$ be a function satisfying
\begin{equation}\label{o23}
  \displaystyle\lim_{n\rightarrow\infty}\sqrt[n]{|\frac{f_{n}}{\Psi_{n}}}|=\tau,\quad \Psi_{n}=\frac{q^{\frac{n(n-1)}{2}}}{[n]_{q}!},\quad0\leq\tau< \min\{\frac{2}{1-q},\frac{j^{(2)}_{1,\alpha}}{1-q}\} ,
\end{equation}
where $j^{(2)}_{1,\alpha}$ is the smallest positive zero of $J_{\alpha}^{(2)}(\cdot;q^{2})$. Then
\begin{equation*}
  f(z)=\sum_{n=0}^{\infty} L_{n}(f)\frac {B^{(2)}_{n,\alpha}(z;q)}{[n]_{q}!},
\end{equation*}
where
\begin{equation}\label{o8}
\begin{split}
   L_{n}(f)=(1-q)^{-n}\sum_{k=n}^{\infty}f_{k}\,\frac{q^{\frac {-k(k-1)}{2}}(q;q)_{k}
(q^{2\alpha+1};q^{2})_{k-n}(\frac{1}{2})^{k-n}}{(q;q)_{k-n}(q^{2\alpha+1};q)_{k-n}}\quad (n\in\mathbb{N}_{0}). \end{split}
\end{equation}

\end{theorem}
\begin{proof}
The generating function of $(B^{(2)}_{n,\alpha}(z;q))_{n=0}^{\infty}$ can be written as
 \begin{equation*}
  B(\xi) E_{q}(z\xi)=\sum_{n=0}^{\infty}B^{(2)}_{n,\alpha}(z;q)\frac{\xi^{n}}{[n]_{q}!},\quad|\xi|<\frac {j^{(2)}_{1,\alpha}}{1-q},
 \end{equation*}
where
 \begin{equation*}
  B(\xi)=\frac {E_{q}(\frac{-\xi}{2})}{g^{(2)}_{\alpha}(i\xi;q)}.
 \end{equation*}
Hence $\Omega=\Delta=\{\xi\in\mathbb{C}:|\xi|<\frac {j^{(2)}_{1,\alpha}}{1-q}\}$.  Therefore, from Theorem \ref{jj},
\begin{equation*}
  f(z)=\sum_{n=0}^{\infty} L_{n}(f)\frac{B^{(2)}_{n,\alpha}(z;q)}{[n]_{q}!},
\end{equation*}
where
\begin{equation}\label{f3}
   L_{n}(f)=\frac{1}{2\pi i}\int_{\Gamma}\frac{w^{n}}{B(w)}F(w)dw,
\end{equation}
$\Gamma$ is a circle $|w|=\rho$ with $\tau<\rho<\min\{\frac{2}{1-q},\frac{j^{(2)}_{1,\alpha}}{1-q}\}$ and  $F(w)$ is defined in (\ref{r32}). From (\ref{r32})  and (\ref{o23}),
\begin{equation}\label{f573}
   L_{n}(f)=\frac{1}{2\pi i}\int_{\Gamma} \left(\sum_{k=0}^{\infty}\frac{f_{k}}{\Psi_{k}w^{k+1}}\right)\frac{w^{n}}{B(w)}dw.
\end{equation}
Since $\dfrac{F(w)}{B(w)}$ is analytic  for $|w|=\rho$, $\tau<\rho<\min\{\frac{2}{1-q},\frac{j^{(2)}_{1,\alpha}}{1-q}\}$, then we can interchange the summation with the integral in (\ref{f573}) to obtain
\begin{equation*}
\begin{split}
   L_{n}(f)&=\sum_{k=n}^{\infty}\frac{f_{k}}{\Psi_{k}}\frac{1}{2\pi i}\int_{\Gamma}\frac{(1/B(w))}{w^{k-n+1}}dw\\&=
   \sum_{k=n}^{\infty}\frac{f_{k}}{\Psi_{k}}\frac{1}{(k-n)!}\frac{d^{k-n}}{dw^{k-n}}\left(\frac{1}{B(w)}\right)|_{w=0}\\&=
   \sum_{k=n}^{\infty}\frac{f_{k}}{\Psi_{k}}\frac{1}{(k-n)!}\frac{d^{k-n}}{dw^{k-n}}\left(g^{(2)}_{\alpha}(iw;q) e_{q}(w/2)\right)|_{w=0},
   \end{split}
\end{equation*}
where we used the Cauchy integral formula, see \cite{cauchy}. Since  $\dfrac{f^{(r)}(0)}{r!}=\dfrac{D^{r}_{q}f(0)}{[r]_{q}!}$, then
\begin{equation}\label{f17}
\begin{split}
   L_{n}(f)= \sum_{k=n}^{\infty}\frac{f_{k}}{\Psi_{k}}\frac{1}{[k-n]_{q}!}D_{q}^{k-n}\,\left(g^{(2)}_{\alpha}(iw;q) e_{q}(w/2)\right)|_{w=0}.\end{split}
\end{equation}
From (\ref{d8}), we conclude that
\begin{equation}\label{f97}
\begin{split}
\frac {D_{q}^{k-n}\left(g^{(2)}_{\alpha}(iw;q)
e_{q}(w/2)\right)|_{w=0}}{[k-n]_{q}!}=\frac{(1-q)^{k-n}(1/2)^{k-n}(q^{\alpha+\frac{1}{2}};q)_{k-n}(-q^{\alpha+\frac{1}{2}};q)_{k-n}}{(q;q)_{k-n}(q^{2\alpha+1};q)_{k-n}}.\end{split}
\end{equation}
Substituting from (\ref{f97}) into (\ref{f17}), we obtain
\begin{equation}\label{g107}
\begin{split}
   L_{n}(f)= \sum_{k=n}^{\infty}\frac{f_{k}}{\Psi_{k}}\frac{(1-q)^{k-n}(1/2)^{k-n}(q^{\alpha+\frac{1}{2}};q)_{k-n}(-q^{\alpha+\frac{1}{2}};q)_{k-n}}
   {(q;q)_{k-n}(q^{2\alpha+1};q)_{k-n}}.\end{split}
\end{equation}
Now, substituting with  $\Psi_{k}=\dfrac{q^{\frac{k(k-1)}{2}}}{[k]_{q}!}$  from (\ref{o23}) into (\ref{g107}), we obtain (\ref{o8}) and complete the proof of this theorem.
\end{proof}
\begin{proposition}\label{sm}
 Let $ f(z)=\displaystyle\sum_{n=0}^{\infty}f_{n}z^{n}$. Assume that the condition of Theorem \ref{bb} holds.\\
(i) If $f$ is a function of $q^{-1}$-exponential growth of order less than one, then $\tau =0$.\\
(ii) If $f$ is a function of $q^{-1}$-exponential growth of order one and type $\gamma$, then $\tau<\frac{q^{\frac{1}{2}-\gamma}}{1-q}$.
\end{proposition}
\begin{proof}
From (\ref{r65}), if $f$ is a function of $q^{-1}$-exponential growth of order $k$ and type $\gamma$ then
\begin{equation*}
  |f_{n}|\leq K \,q^{\frac{(n-\gamma)^{2}}{2k}}\quad (n\in\mathbb{N}_{0}).
\end{equation*}
Hence
\begin{equation}
 \sqrt[n]{|\frac{f_{n}}{\Psi_{n}}|}\leq\sqrt[n]{K \frac {q^{\frac{n^{2}}{2}(\frac{1}{k}-1)-\frac{n}{2}(2\gamma-1)+\frac{\gamma^{2}}{2k}}(q;q)_{n}}{(1-q)^{n}}}\quad (n\in\mathbb{N}_{0}).
\end{equation}
Taking the limit as $n\rightarrow\infty$ yields $\tau=0$ if $f$ is a function of $q^{-1}$-exponential growth of order $k<1$ and $\tau<\frac{q^{\frac{1}{2}-\gamma}}{1-q}$ if $f$ is a function of $q^{-1}$-exponential growth of order one.
\end{proof}
\begin{corollary}\label{oo}
Let  $f(z)=\displaystyle\sum_{n=0}^{\infty}f_{n}z^{n}$ be a function satisfying (\ref{o23}). Then
\begin{equation*}
  f(z)=\sum_{n=0}^{\infty} L_{n}(f)\frac{B_{n}(z;q)}{[n]_{q}!},
\end{equation*}
where $(B_{n}(z;q))_{n}$  defined in (\ref{q8}) and
\begin{equation}\label{o2}
\begin{split}
   L_{n}(f)=(1-q)^{-n}\sum_{k=n}^{\infty}f_{k}\,\frac{q^{\frac {-k(k-1)}{2}}(q;q)_{k}
(-q;q)_{k-n}(\frac{1}{2})^{k-n}}{(q^{2};q)_{k-n}}.
    \end{split}
\end{equation}

\end{corollary}
\begin{proof}
The proof follows  by the substitution $\alpha=\dfrac{1}{2}$ in Theorem \ref{bb}.
\end{proof}
\begin{corollary}\label{oU}
Let  $f(z)=\displaystyle\sum_{n=0}^{\infty}f_{n}z^{n}$ be a function satisfying (\ref{o23}). Then
\begin{equation*}
  f(z)=\sum_{n=0}^{\infty} L_{n}(f)\frac{E_{n}(z;q)}{[n]_{q}!},
\end{equation*}
where $(E_{n}(z;q))_{n}$  defined in (\ref{q13248}) and
\begin{equation*}
\begin{split}
   L_{n}(f)=(1-q)^{-n}\sum_{k=n}^{\infty}f_{k}\,\frac{q^{\frac {-k(k-1)}{2}}(q;q)_{k}
(-1;q)_{k-n}(\frac{1}{2})^{k-n}}{(q;q)_{k-n}}.
    \end{split}
\end{equation*}

\end{corollary}
\begin{proof}
The proof follows from Theorem \ref{bb} by substituting with  $\alpha=\dfrac{-1}{2}$.
\end{proof}
\vskip0.2cm
 Now, we give some examples on Theorem \ref{bb}.
\vskip0.3cm
{\bf{Example 4.1}}
\vskip0.3cm
 If $f(z)=z^{n}$ in Theorem \ref{bb}, then we have

\begin{equation*}
\begin{split}
   L_{m}(f)=\left\{
              \begin{array}{ll}
                q^{\frac {-n(n-1)}{2}}\dfrac{(1-q)^{-m}(q;q)_{n}
(q^{2\alpha+1};q^{2})_{n-m}(\frac{1}{2})^{n-m}}{(q;q)_{n-m}(q^{2\alpha+1};q)_{n-m}}, & \hbox{$ m=0,1,\cdots,n$;} \\
                0, & \hbox{$m > n$.}
              \end{array}
            \right.
 \end{split}
\end{equation*}
 Hence

\begin{equation*}
\begin{split}
z^{n}&=\sum_{k=0}^{n} L_{n-k}(f)\frac{B^{(2)}_{n-k,\alpha}(z;q)}{[n-k]_{q}!}\\&=q^{\frac {-n(n-1)}{2}}\sum_{k=0}^{n} \frac{(1-q)^{k-n}(q;q)_{n}
(q^{2\alpha+1};q^{2})_{k}(\frac{1}{2})^{k}}{(q;q)_{k}(q^{2\alpha+1};q)_{k}}\frac{B^{(2)}_{n-k,\alpha}(z;q)}{[n-k]_{q}!}\\&=q^{\frac {-n(n-1)}{2}}\sum_{k=0}^{n}\left[
                             \begin{array}{c}
                               n \\
                               k \\
                             \end{array}
                           \right]_{q}
\frac{(q^{2\alpha+1};q^{2})_{k}}{2^{k}(q^{2\alpha+1};q)_{k}}B^{(2)}_{n-k,\alpha}(z;q),\end{split}
\end{equation*}
which coincide with the result which we proved in \cite[Theorem. 2.17]{sahar}.
\vskip0.2cm
{\bf{Example 4.2}}
\vskip0.2cm
 Let $f(z)=(z;q)_{n}$. From \cite[P. 25 ]{Gasper},
$$(z;q)_{n}=\displaystyle\sum_{m=0}^{n}\left[
                                           \begin{array}{c}
                                             n \\
                                             m \\
                                           \end{array}
                                         \right]_{q}q^{\frac{m(m-1)}{2}}(-z)^{m}.
$$
Hence by applying Theorem \ref{bb}, we obtain
\begin{equation*}
 (z;q)_{n}= \sum_{m=0}^{n} L_{m}(f)\frac{B^{(2)}_{m,\alpha}(z;q)}{[m]_{q}!},
\end{equation*}
where
\begin{equation}\label{r1879}
\begin{split}
   L_{m}(f)=\left\{
              \begin{array}{ll}
                q^{\frac {-n(n-1)}{2}}\displaystyle\sum_{m=0}^{n}\left[
\begin{array}{c}
                                             n \\
                                             m \\
                                           \end{array}
                                         \right]_{q}\dfrac{(-1)^{m}q^{\frac{m(m-1)}{2}}(q;q)_{n}
(q^{2\alpha+1};q^{2})_{n-m}(\frac{1}{2})^{n-m}}{(1-q)^{m}(q;q)_{n-m}(q^{2\alpha+1};q)_{n-m}}, & \hbox{$0\leq m\leq n$;} \\
                0, & \hbox{otherwise.}
              \end{array}
            \right.
 \end{split}
\end{equation}
 Now,  set $k=n-m$ in (\ref{r1879}) yields
\begin{equation*}
\begin{split}
   L_{n-k}(f)&=(-1)^{n}(1-q)^{-n}\sum_{k=0}^{n}(-1)^{k}\frac{q^{\frac{k^{2}-2nk+k}{2}}(q;q)_{n}(q;q)_{n}
(q^{2\alpha+1};q^{2})_{k}(\frac{(1-q)}{2})^{k}}{(q;q)_{k}(q;q)_{n-k}(q;q)_{k}(q^{2\alpha+1};q)_{k}}\\&=
\frac{(-1)^{n}(q;q)_{n}}{(1-q)^{n}}\sum_{k=0}^{n}\frac{(q^{-n};q)_{k}
(q^{\alpha+\frac{1}{2}};q)_{k}(-q^{\alpha+\frac{1}{2}};q)_{k}(\frac{q(1-q)}{2})^{k}}{(q;q)_{k}(q^{2\alpha+1};q)_{k}(q;q)_{k}}\\&=(-1)^{n}[n]_{q}!\,\,_{3}\phi_{2}(q^{-n},q^{\alpha+\frac{1}{2}},
-q^{\alpha+\frac{1}{2}};q^{2\alpha+1},q;q,\frac{q(1-q)}{2}).
 \end{split}
\end{equation*}
Hence
\begin{equation*}
 (z;q)_{n}= \sum_{k=0}^{n}L_{n-k}(f)\frac{B^{(2)}_{n-k,\alpha}(z;q)}{[n-k]_{q}!},
\end{equation*}
where
\begin{equation*}
\begin{split}
   L_{n-k}(f)=(-1)^{n}[n]_{q}!\,_{3}\phi_{2}(q^{-n},q^{\alpha+\frac{1}{2}},
-q^{\alpha+\frac{1}{2}};q^{2\alpha+1},q;q,\frac{q(1-q)}{2}).
 \end{split}
\end{equation*}


\begin{thebibliography}{10}

\bibitem{cauchy}
L.~Ahlfors.
\newblock {\em {Complex} analysis. {An} introduction to the theory of analytic
  functions of one complex variable}.
\newblock New York-Toronto-London:McGraw-Hill, 1953.

\bibitem{Alsalam1}
W.~A. Al-Salam.
\newblock $ q$-{Appell} polynomials.
\newblock {\em Ann. Mat. Pura Appl.}, 77:31--45, 1967.

\bibitem{Annaby}
M.~H. Annaby and Z.~S. Mansour.
\newblock {\em $q$-Fractional {Calculus} and {Equations}}.
\newblock Lecture Notes in Mathematics 2056. Springer-Verlag., Berlin, 2012.

\bibitem{Appell}
P.~Appell.
\newblock Sur une classe de polyn$\hat{o}$mes.
\newblock {\em Ann. Sci. Ec. Norm. Sup.}, 9:119--144, 1880.

\bibitem{boass1}
R.~P. Boas and R.~C. Buck.
\newblock {\em Entire functions}.
\newblock Academic press, New York, 1954.

\bibitem{boass}
R.~P. Boas and R.~C. Buck.
\newblock {\em {Polynomial} expansions and analytic functions}.
\newblock Springer-Verlag, Berlin, Heidelberg, New York, 1958.

\bibitem{Cardoss}
J.~L. Cardoso.
\newblock Basic {Fourier} series convergence on and outside the $q$-{Linear}
  grid.
\newblock {\em J. Fourier Anal. Appl.}, 17(1):96--114, 2011.

\bibitem{sahar}
S.Z.H. Eweis and Z.S.I. Mansour.
\newblock {Generalized} $q$-{Bernoulli} polynomials generated by {Jackson}
  $q$-{Bessel} functions.
\newblock arxiv.org/abs/2201.10117.

\bibitem{exton}
H.~Exton.
\newblock {\em $q$-{Hypergeometric} {Functions} and {Applications}}.
\newblock {New} {York}, 1983.

\bibitem{Frappier1}
C.~Frappier.
\newblock Representation formulas for entire functions of exponential type and
  generalized {Bernoulli} polynomials.
\newblock {\em J. Austr. Math. Soc. Ser.}, 64:307--316, 1998.

\bibitem{Frappier2}
C.~Frappier.
\newblock Generalized {Bernoulli} polynomials and series.
\newblock {\em Bull. Austral. Math. Soc.}, 61:289--304, 2000.

\bibitem{F3}
C.~Frappier.
\newblock A unified calculus using the generalized {Bernoulli} polynomial.
\newblock {\em J. Approx. Theory.}, 109:279--313, 2001.

\bibitem{Gasper}
G.~Gasper and M.~Rahman.
\newblock {\em Basic {Hypergeometric} {Series}}.
\newblock Cambridge { University} Press, Cambridge, second edition, 2004.

\bibitem{IZ}
M.~E.~H. Ismail and Z.~S. Mansour.
\newblock $q$-analogue of {Lidstone} expansion theorems, two-point {Taylor}
  expansions theorems and {Bernoulli} polynomials.
\newblock {\em Analysis and Apllications.}, 17:853--895, 2019.

\bibitem{Jackson2}
F.~H. Jackson.
\newblock A basic-sine and cosine with with symbolical solutions of certain
  differential equations.
\newblock {\em Proc. Edinb. Math. Soc.}, 22:28--38, 1903/1904.

\bibitem{Jackson}
F.~H. Jackson.
\newblock The basic gamme function and elliptic functions.
\newblock {\em proc. Roy. Soc.}, A 76:127--144, 1905.

\bibitem{Jackson1}
F.~H. Jackson.
\newblock On $q$-functions and certain difference operator.
\newblock {\em Trans. Roy. Soc. Edinb.}, 46:64--72, 1908.

\bibitem{two}
M.~E. Keleshteri and N.~I. Mahmoudov.
\newblock A study on $q$-{Appell} polynomials from determinanatal point of
  view.
\newblock {\em Appl. Math. Comput.}, 260:351--369, 2015.

\bibitem{marim}
Z.~S. Mansour and M.~Al-Towalib.
\newblock The $q$-{Lidstone} series involving $q$-{Bernoulli} and $q$-{Euler}
  polynomials generated by the third {Jackson} $q$-{Bessel} function.
\newblock Submitted.

\bibitem{nac}
L.~Nachbin.
\newblock An extension of the notion of integral functions of finite
  exponential type.
\newblock {\em An. Acad. Brasil. Cienc.}, 16:143--147, 1944.

\bibitem{modified}
M.A. Olshanetsky and V.B.K. Rogov.
\newblock The modified $q$-{Bessel} and the {Bessel}-{Macdonald}{ Functions}.
\newblock {\em arxiv:$q$-alg/9509013.}

\bibitem{olver}
F.~W.~J. Olver.
\newblock {\em {Asymptotics} and special functions}.
\newblock Academic Press, New York., first edition, 1974.

\bibitem{ramis}
J.~P. Ramis.
\newblock About the growth of entire functions solutions of linear algebraic
  $q$-difference equations.
\newblock {\em Ann. Fac. Sci. Toulouse Math.}, 1(16):53--94, 1992.

\bibitem{Sadj}
P.~N. Sadjang.
\newblock On a new q-analogue of {Appell} polynomials.
\newblock {\em Math. CA.}, pages 1--17, 2018.

\end{thebibliography}
\end{document}